# On the logarithmic terms in the asymptotic expansion of integrals[1]

Achim Hennings[2]

**Abstract:** We consider integrals over vanishing cycles in the Milnor fibration of an isolated singularity defined by a Newton non-degenerate function. We single out a condition where the leading logarithmic term of the expansion of the integral into a logarithmic sum can be determined exactly.

## 0    Basic notions and the problem

Let $f:(\mathbb{C}^n, 0) \to (\mathbb{C}, 0)$ be holomorphic with an isolated singularity, $n \geq 2$.

### 0.1    The Milnor fibration and expansion of cohomology classes

We choose a small $\varepsilon > 0$ such that for all $0 < \varepsilon' \leq \varepsilon$ the $\varepsilon'$−sphere $S_{\varepsilon'}$ about $0 \in \mathbb{C}^n$ intersects the fibre $f^{-1}(0)$ transversally. Then we choose $\eta > 0$ small enough, such that $S_\varepsilon$ is still transversal to $f^{-1}(t)$, $|t| \leq \eta$, and the restriction of $f$ to $B_\varepsilon := \{x \in \mathbb{C}^n | \|x\| < \varepsilon\}$ has no critical values $t$, $0 < |t| \leq \eta$.

We introduce the notation $D := D_\eta := \{t \in \mathbb{C} | |t| < \eta\}$, $D^* := D \setminus \{0\}$, $X := B_\varepsilon \cap f^{-1}(D)$, $X^* := X \setminus f^{-1}(0)$, $X_t := f^{-1}(t) \subseteq X$. Then $f: X^* \to D^*$ is a smooth fibre bundle (which does not depend on $\varepsilon$ up to diffeomorphism). $f: X \to D$ is called a Milnor representative.

To the Milnor fibration one associates the local system of the (only interesting degree $n-1$) complex cohomology of the fibres,

$$H^{n-1} := R^{n-1} f_* \mathbb{C}_{X^*} = \bigcup_{t \in D^*} H^{n-1}(X_t, \mathbb{C}).$$

The inverse image sheaf $u^{-1} H^{n-1}$ with respect to the universal covering

$$u: U := \{\tau \in \mathbb{C} | |e^\tau| < \eta\} \to D^*, u(\tau) := e^\tau,$$

is constant. Each fibre is canonically isomorphic to the $\mu$−dimensional vector space $E := \Gamma(U, u^{-1} H^{n-1})$, $\mu = \mu(f, 0)$ being the Milnor number.

On $E$ the monodromy operator $M$ is defined by $MA(\tau) = A(\tau + 2\pi i)$. By primary decomposition, we have $M = M_s M_u = M_u M_s$ where $M_s$ is diagonalisable and $M_u$ is unipotent. There is a unique nilpotent Operator $N$ satisfying $\exp(-2\pi i N) = M_u$.

Let $\Omega^p$ denote the sheaf of holomorphic $p$−forms on $\mathbb{C}^n$. Each $n$−form $\varphi \in H^0(X, \Omega^n)$ gives rise to a holomorphic section of the cohomology bundle by

$$s[\varphi] \in \Gamma(D^*, H^{n-1} \otimes \mathcal{O}_{D^*}), s[\varphi](t) = \left[\frac{\varphi}{df}\Big|X_t\right] \in H^{n-1}(X_t, \mathbb{C}).$$

This section has a convergent development ([Mal])

---

(0-1) $$s[\varphi](t) = \sum_{\lambda \in \Lambda} \sum_{\alpha \in L(\lambda)} \sum_{k=0}^{n-1} A_{\alpha,k} t^\alpha (\ln t)^k / k!,$$

where $\Lambda := $ spectrum $(M)$, $L(\lambda) := \{\alpha > -1 | e^{-2\pi i \alpha} = \lambda\}$, $A_{\alpha,k} \in E$, and this equation is more precisely valid on $U$ (i.e. $\ln t = \tau$, $t = e^\tau$, $A_{\alpha,k} = A_{\alpha,k}(\tau)$). Because $s[\varphi]$ is single-valued on $D^*$, we have ([AGV, p. 354, lemma 13.1])

$$A_{\alpha,k} = N^k A_{\alpha,0} \text{ for all } \alpha, k.$$

Furthermore $N^{n-1} A_{\alpha,0} = 0$ for $\alpha \in \mathbb{Z}$ ([AGV, p. 379, cor. 2]).

The leading contribution in (0-1) obviously is given by the summand $A_{\alpha,k} \neq 0$ with largest $k$ among those with smallest $\alpha$. Therefore, we define $(\alpha, k)(\varphi) = (\alpha(\varphi), k(\varphi))$ with

$$\alpha(\varphi) := \min\{\alpha | \exists k : A_{\alpha,k} \neq 0\}$$

$$k(\varphi) := \max\{k | A_{\alpha(\varphi),k} \neq 0\}.$$

Then $(\alpha, k)(\varphi) = \min\{(\alpha, k) | A_{\alpha,k} \neq 0\}$ for the lexicographic order

$$(\alpha, k) \leq (\alpha', k') \Leftrightarrow \alpha < \alpha' \text{ or } (\alpha = \alpha' \text{ and } k \geq k').$$

### 0.2 The canonical Hermitian form of Barlet

We denote the sheaf of smooth $p$-forms on $\mathbb{C}^n$ by $\mathcal{E}^p$. In [Ba 1], D. Barlet has shown that for $\varphi \in \Gamma(X, \mathcal{E}^{2n})$ and an auxiliary cut-off function $\sigma \in C_c^\infty(X)$, $\sigma = 1$ near $0 \in \mathbb{C}^n$, the real fibre integral

$$F(t) := F_{\sigma\varphi}(t) := \int_{X_t} \frac{\sigma\varphi}{df \wedge d\bar{f}}, t \in D^*,$$

admits an asymptotic expansion[3][4]

$$F(t) \sim F^* := \sum_{\substack{\alpha \in \mathcal{A} \\ m, m' \in \mathbb{N} \\ k = 0, \ldots, n-1}} c_{m,m'}^{\alpha,k} t^m \bar{t}^{m'} |t|^{2\alpha} (\ln|t|^2)^k,$$

where $\mathcal{A} := \{\alpha \in ]-1, 0] | \exp(-2\pi i \alpha) \in \Lambda\}$ and $c_{m,m'}^{\alpha,k} \in \mathbb{C}$.

The residue class $F^* \mod \mathbb{C}[\![t, \bar{t}]\!]$ in the formal vector space

$$\mathcal{N} := \bigoplus_{\substack{\alpha \in ]-1,0] \\ k = 0, \ldots, n-1}} \mathbb{C}[\![t, \bar{t}]\!] |t|^{2\alpha} (\ln|t|^2)^k / \mathbb{C}[\![t, \bar{t}]\!]$$

is independent of $\sigma$, as also $F$ does not depend on $\sigma$ up to a smooth function on $D$.

We obtain a map

$$\mathcal{H}: \Omega_0^n \times \Omega_0^n \to \mathcal{N}$$

---

[3] Meaning that $\forall a: F(t) - \sum_{m+m'+2\alpha \leq a} c_{m,m'}^{\alpha,k} t^m \bar{t}^{m'} |t|^{2\alpha} (\ln|t|^2)^k = O(|t|^{a+\varepsilon})$ for some $\varepsilon > 0$.

[4] At first sight, only $f^n \bar{f}^n \varphi$ is divisible by $df \wedge d\bar{f}$, and $|t|^{2a} F$, $a > n - 1$, has an expansion of the stated kind. But we have also $f^n \varphi = df \wedge \alpha$ with $\alpha$ of type $(n-1, n)$, hence $\alpha = \bar{\partial} \beta$, and

$$\frac{\partial}{\partial \bar{t}} \int_{X_t} \sigma \beta = \int_{X_t} \frac{df \wedge \bar{\partial}(\sigma\beta)}{df \wedge d\bar{f}} \equiv \int_{X_t} \frac{\sigma f^n \varphi}{df \wedge d\bar{f}} \mod C^\infty(D).$$

Because the asymptotic expansion of $\int_{X_t} \sigma\beta$ can be differentiated term by term ([Ba 1]), all terms in $t^n F^*$ and $F^*$ are of order $> -1$ in $\bar{t}$. The same is true for $t$, and the claimed expansion follows.

$$\mathcal{H}(\omega, \omega') = \left(\frac{-1}{2\pi i}\right)^{n-1} F^*_{\sigma\varphi} \bmod \mathbb{C}[\![t, \bar{t}]\!],$$

where $\varphi = \sigma\widetilde{\omega} \wedge \overline{\widetilde{\omega}'}$ and $\widetilde{\omega}, \widetilde{\omega}'$ are holomorphic representatives in a neighborhood of the support of the cut-off function $\sigma$.

In [Ba 2] it is shown that this map can in fact be defined on $L \times L$, $L := \Omega_0^n / df \wedge d\Omega_0^{n-2}$, and then becomes $\mathbb{C}\{t\}$ –Hermitian and quasi-horizontal, i.e. denoting the covariant derivative

$$\nabla_{\partial/\partial t} : L \to L \otimes Q(\mathbb{C}\{t\})$$

by $\nabla_t$, we have

$$\mathcal{H}([\omega], [\omega']) = \overline{\mathcal{H}([\omega'], [\omega])},$$

$$\mathcal{H}(., [\omega']) \text{ is } \mathbb{C}\{t\} \text{ –linear,}$$

$$\mathcal{H}(\nabla_t[\omega], [\omega']) = \frac{\partial}{\partial t}\mathcal{H}([\omega], [\omega']) \text{ if } \nabla_t[\omega] \in L.$$

From these properties, D. Barlet in [Ba 2] infers the existence of an ordinary Hermitian form $h: E \times E \to \mathbb{C}$, which satisfies

(i) $h(MA, MA') = h(A, A')$      (Monodromy invariance)
(ii) $h(NA, A') = h(A, NA')$      (Self-Adjointness of $N$)
    (these two conditions are equivalent, cf. [Ba 2, append.])
(iii) If the sections of $\omega, \omega' \in \Gamma(X, \Omega^n)$ are expanded as
$$s[\omega^{(\prime)}] = \sum_\alpha \sum_k N^k A_\alpha^{(\prime)} t^\alpha (\ln t)^k / k!,$$
then $\mathcal{H}([\omega], [\omega']) \in \mathcal{N}$ is represented as
$$\mathcal{H}([\omega], [\omega']) = \sum_{\alpha, \alpha' \notin \mathbb{Z}, \alpha - \alpha' \in \mathbb{Z}} \sum_{k=0}^{n-1} c_{\alpha, \alpha', k} \, t^\alpha \bar{t}^{\alpha'} (\ln|t|^2)^k / k! +$$
$$\sum_{\alpha, \alpha' \in \mathbb{Z}} \sum_{k=0}^{n-2} c_{\alpha, \alpha', k} \, t^\alpha \bar{t}^{\alpha'} (\ln|t|^2)^{k+1} / (k+1)!,$$
where $c_{\alpha, \alpha', k} = h(N^k A_\alpha, A'_{\alpha'})$ (cf. [Loe]).

Note that $N^{n-1} A_\alpha = 0$ for $\alpha \in \mathbb{Z}$ (as part of the monodromy theorem) and that $c_{\alpha, \alpha', k} = 0$ if $\alpha - \alpha' \notin \mathbb{Z}$, since primary components of $M$ are orthogonal with respect to $h$ (cf. [Ba 2, append.]) .

### 0.3     The Newton polyhedron

The following notation is common:

For $g \in \mathbb{C}[\![x]\!]$, $x = (x_1, \ldots, x_n)$, $g = \sum g_m x^m$, and $A \subseteq \mathbb{R}_+^n$ we write

$$supp\,(g) := \{m \in \mathbb{N}^n | g_m \neq 0\},$$

$$\Gamma_+(g) := conv\,(supp\,(g) + \mathbb{N}^n) \quad \text{(Newton polyhedron)}$$

$$\Gamma(g) := \text{union of all compact faces of } \Gamma_+(g) \quad \text{(Newton boundary)}$$

$$g_A := \sum_{m \in A} g_m x^m,$$

$$L(A) := \{h \in \mathbb{C}[x] | supp\,(h) \subseteq A\}.$$

$g$ is called *non-degenerate* (in the sense of [Kou, 1.19]) if for all compact faces $\delta$ of $\Gamma_+(g)$ the (quasi-homogeneous) polynomial $g_\delta$ has no critical points in $(\mathbb{C}^*)^n$ (equivalently the divisor of $g_\delta$ is nonsingular). In [Kou] $g$ is called convenient (commode) if $\Gamma_+(g)$ meets all coordinate axes. (For isolated singularities this condition usually can be assumed by adding a suitable homogeneous polynomial of high degree.)

We consider for a fixed non-unit power series $g$, which is convenient, the polyhedra $\Gamma_+ := \Gamma_+(g)$, $\Gamma := \Gamma(g)$. Then we put for any nonzero $h \in \mathbb{C}[\![x]\!]$

$$v(h) := \sup\{a \in \mathbb{R}_+ | supp\ (h) \subseteq a\Gamma_+\} \quad \text{(Newton order)}$$

$$l(h) := \max\left\{l \in \mathbb{N} \,\middle|\, \begin{matrix}\text{there is a compact face } \delta \text{ of } \Gamma_+ \text{ of dimension } n-1-l \\ \text{such that } supp\ (h) \cap v(h)\delta \neq \emptyset\end{matrix}\right\}.$$

(To a unit is attached the Newton order 0 by the convention $0 \cdot \Gamma_+ = \mathbb{R}_+^n$.)

For holomorphic $p$ —forms (with at most logarithmic poles along $(x_1 ... x_n)$) we use the representation in the logarithmic basis

$$\varphi = \sum_{|I|=p} \varphi_I (dx)_I / x_I,$$

where for $I = \{i_1 < \cdots < i_p\} \subseteq [1, n]$ we have written $(dx)_I = dx_{i_1} \wedge ... \wedge dx_{i_p}$ and $x_I = x_{i_1} ... x_{i_p}$. Then we define the support as $supp\ (\varphi) := \bigcup_I supp\ (\varphi_I)$ and similarly $\varphi_A := \sum_{|I|=p} \varphi_{I,A} (dx)_I / x_I$. For the definitions of $\Gamma_+(\varphi), \Gamma(\varphi), v(\varphi), l(\varphi)$ we merely have to replace the support of a function by $supp\ (\varphi)$.

## 0.4 Content of this article

From now on we assume that the function $f$ is non-degenerate and convenient with Newton polyhedron $\Delta := \Gamma_+(f)$, and with respect to $\Delta$ we form the pair of numbers $(v(\varphi), l(\varphi))$ for holomorphic $n$ —forms.

By a result of V.A. Vasil'ev [Vas] the lowest contribution $(\alpha(\varphi), k(\varphi))$ to the expansion (0-1) can be estimated in the lexicographical order as

(0-2)  $(\alpha, k)(\varphi) \geq (v, l)(\varphi) - (1,0)$ for all $\varphi \in \Gamma(X, \Omega^n)$.

*In this article, we establish a sufficient condition on $\varphi$ for equality in* (0-2). The main result is Th. 3.9.

We consider a compact face $\delta$ of $\Delta$ of dimension $n - r$, $1 \leq r \leq n$, (which is not contained in a coordinate hyperplane) and a certain positive multiple $a\delta$. Then we consider at the same time two polynomial $n$ —forms in the interior of $a\delta$, $\varphi_1$ taking on the role of $\varphi$,

$$\varphi_j \in \Omega_0^n \text{ with } supp(\varphi_j) \subseteq a\delta^\circ, j = 1,2.$$

As a first case, we suppose that $a$ is integral: $a \in \mathbb{N}^+$ and $1 \leq r \leq n - 1$. The second case, $a \notin \mathbb{Z}$ and $1 \leq r \leq n$, will later be reduced to the first. Obviously, $(v, l)(\varphi_j) = (a, r - 1)$, if $\varphi_j \neq 0$.

We take a cut-off-function $\sigma \in C_c^\infty(X)$, $\sigma = 1$ near $0 \in \mathbb{C}^n$, and form the integral

$$M(\lambda) := M[\sigma\varphi_1 \wedge \bar{\varphi}_2](\lambda) := \int_X |f|^{2\lambda} \sigma\varphi_1 \wedge \bar{\varphi}_2.$$

$M(\lambda)$ is defined and holomorphic for $\lambda \in \mathbb{C}$, $Re\,\lambda > 0$ (as $|f|^{2\lambda} \ln|f|^2$ is continuous). In section 2.1 (together with section 1) we justify that $M(\lambda)$ has a meromorphic continuation to $\mathbb{C}$.

The nature of the poles of $M(\lambda)$ can be determined from the asymptotic expansion $\mathcal{H}(\varphi_1, \varphi_2)$ of the fibre integral $F_{\sigma\varphi_1 \wedge \bar{\varphi}_2}$ by means of Fubini's theorem and the Mellin transform:

Since $X$ is Stein, $\varphi_j/df$ can be viewed as a (non-unique) meromorphic $n-$form, which becomes holomorphic by multiplication with $f^N$ for some $N$ (e.g. $N = n$).

For $Re\,\lambda \geq N$ we have

$$\int_X |f|^{2\lambda} \sigma\varphi_1 \wedge \bar{\varphi}_2 = (-1)^{n-1} \int_X |f|^{2\lambda} \sigma df \wedge d\bar{f} \wedge \frac{\varphi_1}{df} \wedge \frac{\overline{\varphi_2}}{df}$$

$$= (-1)^{n-1} \int_{\mathbb{C}} \left( \int_{X_t} |f|^{2\lambda} \sigma \frac{\varphi_1}{df} \wedge \frac{\overline{\varphi_2}}{df} \right) dt \wedge d\bar{t}$$

$$= (-1)^n \int_0^\infty r^{2(\lambda+1)} \int_{|t|=r} \left( \int_{X_t} \sigma \frac{\varphi_1}{df} \wedge \frac{\overline{\varphi_2}}{df} \right) \frac{dt}{t} \frac{2dr}{r} = (-1)^n \mathcal{M}\phi(\lambda+1)$$

where

$$\phi(r) := \int_{|t|=r} \left( \int_{X_t} \sigma \frac{\varphi_1}{df} \wedge \frac{\overline{\varphi_2}}{df} \right) \frac{dt}{t}$$

and

$$\mathcal{M}\phi(\lambda) := \int_0^\infty r^{2\lambda} \phi(r) \frac{2dr}{r}$$

is the Mellin transform.

The asymptotic expansion of $\phi(r)$ is obtained by integration of each term of $\mathcal{H}(\varphi_1, \varphi_2)$:

$$\phi(r) \sim (2\pi i)^n \left( \sum_{\alpha \notin \mathbb{Z}} \sum_{k=0}^{n-1} c_{\alpha,\alpha,k} r^{2\alpha} \frac{(\ln r^2)^k}{k!} + \sum_{\alpha \in \mathbb{Z}} \sum_{k=0}^{n-2} c_{\alpha,\alpha,k} r^{2\alpha} \frac{(\ln r^2)^{k+1}}{(k+1)!} \right) mod\ \mathbb{C}[\![r^2]\!].$$

This expansion determines the poles of the Mellin transform of $\phi$ by the following lemma.

**Lemma 0.1:** Let $\phi \in C^\infty(\mathbb{R}^+)$ have bounded support and the asymptotic expansion

$$\phi(r) \sim \sum_{\alpha \in A} \sum_{k=0}^K c_{\alpha,k} r^{2\alpha} \frac{(\ln r^2)^k}{k!} \quad (r \to 0),$$

with a subset $A \subseteq \mathbb{R}$, which is discrete closed and bounded below. Then the Mellin transform

$$\mathcal{M}\phi(\lambda) := \int_0^\infty r^{2\lambda}\phi(r)\frac{2dr}{r} \quad (Re\,\lambda > -\min A)$$

has a meromorphic continuation to $\mathbb{C}$ with pole set $\subseteq -A$ and the principal part

$$\sum_{k=0}^{K}(-1)^k c_{\alpha,k}\left(\frac{1}{\lambda+\alpha}\right)^{k+1} \quad \text{at} -\alpha \text{ for } \alpha \in A.$$

Proof: This is an application of the formula

$$\int_0^1 r^{2\lambda}r^{2\alpha}\frac{(\ln r^2)^k}{k!}\frac{2dr}{r} = (-1)^k\left(\frac{1}{\lambda+\alpha}\right)^{k+1} \quad (Re\,\lambda > -\alpha),$$

which is obtained from the simple case $k=0$ by taking derivatives in the variable $\lambda$.

We can conclude that the function $M$ has a meromorphic continuation to $\mathbb{C}$ with rational poles of order $\leq n$, and the principal parts, except the residues at $\alpha \in -1-\mathbb{N}$, are determined by $s[\varphi_j]$ (and independent from $\sigma$). The non-trivial poles $\alpha$ are those with $\alpha \notin \mathbb{Z}$ or pole order $\geq 2$. Let $P_m(\alpha,M)$ denote the Laurent coefficient of $M$ of order $m$ at $\alpha \in \mathbb{C}$.

By the estimate (0-2) the largest possible non-trivial pole is $-a$, and the leading coefficient is

$$P_{-(r-1)}(-a,M) = (-1)^{n+r}(2\pi i)^n h(A^{\varphi_1}_{a-1,r-1}, A^{\varphi_2}_{a-1,0}),$$

provided this is non-zero. Then in particular $A^{\varphi_1}_{a-1,r-1} \neq 0$ (and $A^{\varphi_2}_{a-1,0} \neq 0$).

If we can find $\varphi_2$ such that $P_{-(r-1)}(-a,M) \neq 0$, the estimate (0-2) holds with equality. In section 2, we determine for this purpose $P_{-(r-1)}(-a,M)$ in a second independent way by residue calculations on the defining integral for $M$. This leads to the criterion Cor. 3.6 in section 3.

So far $a$ was assumed integral. Fortunately, the non-integral case can be reduced to the former by a method due to A.N. Varchenko ([AGV]): One replaces $f$ by $f + y^e$ with a new variable $y$ and suitable integer $e$. Also $\varphi$ is multiplied by a power of $y$ in order to receive an integer degree.

## 1      Resolution of $f$

We keep the definitions of section 0: The function $f$ is non-degenerate and convenient with Newton polyhedron $\Delta$, $\delta$ is a compact face of dimension $n-r$ ($1 \leq r \leq n-1$), which is not contained in a coordinate hyperplane, and $\varphi_1, \varphi_2 \in \Omega_0^n$ are polynomial forms with $supp(\varphi_1) \subseteq a\delta°$, $supp(\varphi_2) \subseteq a\delta$, where $a \in \mathbb{N}^+$. (We have weakened here the condition $supp(\varphi_2) \subseteq a\delta°$ only for this section for better explanation.)

We construct a toric resolution of $f$ in the familiar way (cf. [AGV]), using notation as in [Da 1]. Let $\Sigma_0$ be the fan on $\mathbb{R}_+^n$ which is dual to $\Delta$. It consists of the dual cones $\check{\sigma}(\gamma)$ of the cones $\sigma(\gamma) := \mathbb{R}_+(\Delta - p)$, where $\gamma$ is a face of $\Delta$ (including $\Delta$ itself) and $p$ an arbitrary point of $\gamma°$.

Let $\Sigma$ be an arbitrary regular subdivision, which does not subdivide the boundary cones of $\mathbb{R}_+^n$.

We need the set $S$ of cones $\check{\tau} \in \Sigma$, which are contained in $\check{\sigma}(\delta)$ and have the same dimension, so that $\check{\sigma}(\delta) = \bigcup_{\check{\tau} \in S} \check{\tau}$. We obtain a diagram of corresponding torus embeddings:

(1-1)
$$\begin{array}{ccccc} \mathbb{C}^n & \xleftarrow{\rho_1} & X_{\Sigma_0} & \xleftarrow{\rho_2} & X_\Sigma \\ & & \cup| & & \cup| \\ & & X_{\sigma(\delta)} & \leftarrow & \bigcup_{\check{\tau} \in S} X_\tau \end{array}$$

$\bigcup_{\check{\tau} \in S} X_\tau$ is the total inverse image $\rho_2^{-1}(X_{\sigma(\delta)})$. By the choice of $\Sigma$, $D := (\rho_1 \rho_2)^{-1}(0)$ is the exceptional set, consisting of all divisors corresponding to inner rays in $\Sigma$. The strict transforms $Z_0$, resp. $Z$, of the divisor $(f) \subseteq X$ in $X_{\Sigma_0}$, resp. $X_\Sigma$, (restricted to the inverse image of $X$) are hypersurfaces. $Z \cup D$ is a divisor with simple normal crossings.

The affine toric varieties $X_{\sigma(\delta)}$ and $X_\tau$ contain the closed strata (= orbits of the torus-action) $T_\delta := X_{cospan(\sigma(\delta))}$ and $T_\tau := X_{cospan(\tau)}$, which are compatible with $\rho_2$:

(1-2)
$$\begin{array}{ccc} X_{\sigma(\delta)} & \leftarrow & \bigcup_{\check{\tau} \in S} X_\tau \\ \cup| & & \cup| \\ T_\delta & \leftarrow & \bigcup_{\check{\tau} \in S} T_\tau \end{array}$$

**Lemma 1.1:** $\Pi_\tau := \rho_2|_{T_\tau} : T_\tau \to T_\delta$ is an isomorphism.

Proof: The diagram corresponds to the algebraic diagram:

$$\begin{array}{ccc} \mathbb{C}[\sigma(\delta) \cap \mathbb{Z}^n] & \hookrightarrow & \mathbb{C}[\tau \cap \mathbb{Z}^n] \\ \downarrow & & \downarrow \\ \mathbb{C}[cospan(\sigma(\delta)) \cap \mathbb{Z}^n] & = & \mathbb{C}[cospan(\tau) \cap \mathbb{Z}^n] \end{array}$$

The vertical maps are the restrictions defined by $x^m \mapsto x^{\chi(m)m}$, where $\chi$ is the characteristic function of $cospan(\sigma(\delta)) = cospan(\tau)$.

**Lemma 1.2:** The closures $\overline{T}_\tau$ in $X_\Sigma$ are pairwise disjoint.

Proof: $\overline{T}_\tau$ is the intersection of all divisors $D(l) \subseteq X_\Sigma$, which belong to the edges of $\check{\tau}$. For $\check{\tau}_1, \check{\tau}_2 \in S$ the following is true: The intersection of all the divisors $D(l)$ belonging to edges of $\check{\tau}_1$ or $\check{\tau}_2$ is non-empty $\Leftrightarrow$ there is a cone $\sigma \in X_\Sigma$, which possesses all these edges, i.e. contains $\check{\tau}_1, \check{\tau}_2$. This cannot happen, since $\check{\tau}_1, \check{\tau}_2$ both span the plane $\mathbb{R}\check{\sigma}(\delta)$ of dimension $r$.

We write $X'_\Sigma$ for the inverse image of $X$ in $X_\Sigma$. Let $U_\tau \subseteq X_\Sigma$ be disjoint open neighborhoods of $\overline{T}_\tau$, $\check{\tau} \in S$. As $\overline{T}_\tau$ is exceptional, we may choose $U_\tau \subseteq X'_\Sigma$. As above $a \in \mathbb{N}^+$.

Then we have (denoting pulled back differential forms by the same symbol)

$$f^{-a}\varphi_i \in \Gamma(X'_\Sigma, \Omega^n(\log D)(aZ)), \; i = 1,2.$$

Consider one of the $\check{\tau} \in S$. We may form the iterated Poincaré residue of $f^{-a}\varphi_i$ along all $D(l)$, $l \subseteq \check{\tau}$ an edge, (in any fixed order) to obtain

$$Res_{T_\tau}(f^{-a}\varphi_i) \in \Gamma(T_\tau, \Omega^{n-r}(aZ \cap T_\tau)).$$

For $i = 1$ this form is holomorphic along $\overline{T}_\tau \setminus T_\tau$, for $i = 2$ it has at most logarithmic poles.

Now we fix $\check{\tau}_0 \in S$ for reference. For $\check{\tau} \in S$ we have the maps of Lemma 1.1

$$\Pi_{\tau_0}: T_{\tau_0} \xrightarrow{\sim} T_\delta \xleftarrow{\sim} T_\tau : \Pi_\tau,$$

which restrict to $T_\tau \cap Z \xrightarrow{\sim} T_\delta \cap Z_0$. The composition $\Pi_\tau^{\tau_0} := \Pi_\tau^{-1} \Pi_{\tau_0}$ restricts to $T_{\tau_0} \setminus Z \xrightarrow{\sim} T_\tau \setminus Z$.

**Lemma 1.3:** For $\check{\tau} \in S$ the residues

$$\left(\Pi_\tau^{\tau_0}\right)^* Res_{T_\tau}(f^{-a}\varphi_i) \in \Gamma(T_{\tau_0}, \Omega^{n-r}(aZ \cap T_{\tau_0}))$$

are different only by constant factors $\in \mathbb{Q}^*$.

Proof: Let $m \in \delta \cap \mathbb{Z}^n$, $\tilde{f} := x^{-m}f$, $\tilde{\varphi}_i := x^{-ma}\varphi_i$, so that $f^{-a}\varphi_i = \tilde{f}^{-a}\tilde{\varphi}_i \in \tilde{f}^{-a}\mathbb{C}[\tau \cap \mathbb{Z}^n] \otimes \Lambda^n \mathbb{C}^n$. The residue $\tilde{f}^a|_{T_\tau} Res_{T_\tau}(f^{-a}\varphi_i) \in \mathbb{C}[cospan(\tau) \cap \mathbb{Z}^n] \otimes \Lambda^{n-r}\mathbb{C}^n$ is represented as the iterated contraction with the vector fields ($\in \mathbb{C}^{n\vee}$) corresponding to primitive edge vectors of $\check{\tau}$ on the second factor combined with the restriction map on the first factor. In this representation, $\left(\Pi_\tau^{\tau_0}\right)^*$ is the identity. Now we have only to apply the following remark to the edge vectors of different $\check{\tau}$.

**Lemma 1.4:** Let $v_1, \ldots, v_r \in \mathbb{Q}^{n\vee}$ and $w_1, \ldots, w_r \in \mathbb{Q}^{n\vee}$ be linear independent linear forms on $\mathbb{Q}^n$. They define contractions on $\Lambda^\bullet \mathbb{Q}^n$.

(1) $\bigcap_{j=1}^r \ker v_j = \bigcap_{j=1}^r \ker w_j \Leftrightarrow \sum_{j=1}^r \mathbb{Q}v_j = \sum_{j=1}^r \mathbb{Q}w_j$
(2) $w_k = \sum_{j=1}^r a_{jk} v_j$ ($k = 1, \ldots, r$) $\Rightarrow i_{w_1} \circ \ldots \circ i_{w_r} = \det(a_{jk}) i_{v_1} \circ \ldots \circ i_{v_r}$

## 2 Meromorphic continuation of integrals of Mellin type

### 2.1 Local consideration

Let $W \subseteq \mathbb{C}^l$ be an open set. We consider a smooth function $\kappa \in C^\infty(\mathbb{C}^n \times W)$, which is holomorphic in the second variable $w \in W$ and has a proper support for the projection $pr_2: \mathbb{C}^n \times W \to W$. We are concerned with the integral

(2-1) $\quad M(\lambda, w) := \int_{\mathbb{C}^n} |z_1|^{2\lambda_1} \ldots |z_n|^{2\lambda_n} \kappa(z, w) dV$

(with $dV := dx_1 \wedge dy_1 \wedge \ldots \wedge dx_n \wedge dy_n$).

**Lemma 2.1:**

(1) $M(\lambda, w)$ is holomorphic in $(\lambda, w)$ for $Re \, \lambda_i > -1$, $1 \leq i \leq n$.
(2) For $k \leq n$ and $a_1, \ldots, a_k \in \mathbb{N}^+$ we have

$$M(\lambda,w)(-1)^k \prod_{i=1}^{k}[(\lambda_i+1)^2 \ldots (\lambda_i+a_i-1)^2(\lambda_i+a_i)] =$$
$$\int_{\mathbb{C}^n}\left[\prod_{i=1}^{k}|z_i|^{2(\lambda_i+a_i)}\right]\left[\prod_{i=k+1}^{n}|z_i|^{2\lambda_i}\right]\frac{1}{z_1\ldots z_k}\frac{\partial^k}{\partial\bar{z}_1\ldots\partial\bar{z}_k}\prod_{i=1}^{k}\left[\frac{\partial^2}{\partial z_i\partial\bar{z}_i}\right]^{a_i-1}\kappa(z,w)dV.$$

Proof: (1): By assumption, we are allowed to differentiate under the integral sign. The derivatives $\partial/\partial\bar{\lambda}_i$, $\partial/\partial\bar{w}_j$ are zero.

(2): Observe that (for $z, \lambda \in \mathbb{C}$, $\text{Re } \lambda > 0$)

$$\frac{\partial}{\partial z}(z|z|^{2\lambda}) = (\lambda+1)|z|^{2\lambda}, \frac{\partial}{\partial\bar{z}}(\bar{z}|z|^{2\lambda}) = (\lambda+1)|z|^{2\lambda}.$$

Using this, we get by partial integration

$$-M(\lambda,w)(\lambda_i+1) = \int_{\mathbb{C}^n}|z_1|^{2\lambda_1}\ldots|z_i|^{2\lambda_i}\ldots|z_n|^{2\lambda_n}z_i\frac{\partial}{\partial z_i}\kappa(z,w)dV,$$

$$-M(\lambda,w)(\lambda_i+1) = \int_{\mathbb{C}^n}|z_1|^{2\lambda_1}\ldots|z_i|^{2\lambda_i}\ldots|z_n|^{2\lambda_n}\bar{z}_i\frac{\partial}{\partial\bar{z}_i}\kappa(z,w)dV.$$

By repeated application of these formulas, we prove (2).

The right-hand side in (2) is holomorphic for $\text{Re } \lambda_i + a_i > -\frac{1}{2}$, $1 \leq i \leq k$, $\text{Re } \lambda_i > -1$, $k+1 \leq i \leq n$. If we set $\lambda_i + a_i = 0$, $1 \leq i \leq k$, we can apply Cauchy's integral formula for smooth functions to obtain

$$(-\pi)^k \int_{\{z|z_i=0, 1\leq i\leq k\}}|z_{k+1}|^{2\lambda_{k+1}}\ldots|z_n|^{2\lambda_n}\prod_{i=1}^{k}\left[\frac{\partial^2}{\partial z_i\partial\bar{z}_i}\right]^{a_i-1}\kappa(z,w)dV(z_{k+1},\ldots,z_n).$$

Thus, we have proved:

**Corollary 2.2:**

(1) $M(\lambda,w)$ is meromorphic (extendable) on $\mathbb{C}^n \times W$ with at most simple poles at the hyperplanes $\lambda_i = -r$, $r \in \mathbb{N}^+$.
(2) For $\text{Re } \lambda_i > -1$, $k+1 \leq i \leq n$, the coefficient of $\prod_{i=1}^{k}\frac{1}{\lambda_i+a_i}$ in the Laurent-series of $M(\lambda,w)$ at $(-a_1,\ldots,-a_k,\lambda_{k+1},\ldots,\lambda_n,w)$ equals
$$\pi^k\prod_{i=1}^{k}\left[\frac{1}{(a_i-1)!}\right]^2\int_{\mathbb{C}^{n-k}}|z_{k+1}|^{2\lambda_{k+1}}\ldots|z_n|^{2\lambda_n}\prod_{i=1}^{k}\left[\frac{\partial^2}{\partial z_i\partial\bar{z}_i}\right]^{a_i-1}\kappa(z,w)dV(z_{k+1},\ldots,z_n).$$
It depends holomorphically on $(\lambda_{k+1},\ldots,\lambda_n,w)$.

**Remark 2.3:** If the function $\kappa$ is zero on $\{z|z_1\ldots z_k = 0\}$, then $M(\lambda,w)$ has no pole along the hyperplanes $\lambda_i = -1$, $1 \leq i \leq k$.

Proof: By the mean value theorem $|\kappa| \leq \max_u\left(\left|\frac{\partial\kappa}{\partial z_i}(u,w)\right| + \left|\frac{\partial\kappa}{\partial\bar{z}_i}(u,w)\right|\right)|z_i|$ for $1 \leq i \leq k$. Therefore $M(\lambda,w)$ is holomorphic for $\text{Re } \lambda_i > -\frac{3}{2}$, $1 \leq i \leq k$, $\text{Re } \lambda_j > -1$, $k+1 \leq j \leq n$.

## 2.2 Globalization

Let $Y$ be a complex manifold (with countable basis, hence paracompact) and $D$ an effective divisor on $Y$.

By a *tube function for D* we understand a smooth function $g: Y \to \mathbb{R}_+$ which is locally representable as $g = A|t|^2$ with a smooth positive function $A$ and a local equation $t$ for $D$.

If we have a (locally finite) sum $D = \sum_{i \in I} v_i D_i$ with $v_i \geq 0$, we may also speak of a tube function for $\bigcup_{i \in I} D_i$ with multiplicities $v_i$ along $D_i$.

We use the notation $D(J) := \bigcup_{i \in J} D_i$, $D_J := \bigcap_{i \in J} D_i$, $J \subseteq I$.

### 2.2.1 Favorite Situation I

We investigate a situation, which arises when we have an embedded resolution of a hypersurface and we are given (semi-) meromorphic differential forms with poles along the strict transform and certain exceptional divisors. (See appendix 1 for the properties of semi-meromorphic forms.)

With $Y$ a manifold of dimension $n$, let $Z = D_0, D_1, \ldots, D_K$ be smooth divisors on $Y$ with normal crossings, and assume $D_0 \cap \ldots \cap D_k \neq \emptyset$ for some $k \in [1, K]$. Let $I := [0, K]$, $J = [1, k]$. Suppose, we are given a tube function $g$ for $D(I) \subseteq Y$ with multiplicities $v_i > 0$, $i \in [0, k]$, $v := \text{lcm}\{v_0, \ldots, v_k\}$, and $v_i \geq 0$, $i \in [k+1, K]$. For $D_0 \cap D_J \subseteq D_J$ we are given a tube function $\tilde{g}$ with multiplicity $\tilde{v} > 0$.

If we have semi-meromorphic differential forms $\psi_j \in \Gamma(Y, \mathcal{E}_Y^n(\log D(J))(*Z))$, $j = 1,2$, and a smooth function with compact support $\sigma \in C_c^\infty(Y)$ we can form the integrals of Mellin type:

$$M(\lambda) := \int_Y g^\lambda \sigma \psi_1 \wedge \overline{\psi_2},$$

$$\widetilde{M}(\lambda) := \int_{D_J} \tilde{g}^\lambda \sigma \text{Res}_{D_J} \psi_1 \wedge \overline{\text{Res}_{D_J} \psi_2}.$$

The notation $\text{Res}_{D_J} \psi_j$ means the iterated Poincaré residue along the divisors $D_1, \ldots, D_k$ with some fixed order. Since the term appears twice, there is no ambiguity.

**Lemma 2.4:**

(1) $M(\lambda)$ is meromorphic on $\mathbb{C}$ (with poles $\in \frac{1}{v}\mathbb{Z}$ of order $\leq n$).
(2) $\lambda^{k+1} M(\lambda)$ is holomorphic at 0.
(3) $\widetilde{M}(\lambda)$ is meromorphic on $\mathbb{C}$ (with poles $\in \frac{1}{\tilde{v}}\mathbb{Z}$ of order $\leq 1$).
(4) $[\lambda^{k+1} M(\lambda)]_{\lambda=0} = \tilde{v}(v_0 \ldots v_k)^{-1}(2\pi i)^k (-1)^{kn-k(k-1)/2} [\lambda \widetilde{M}(\lambda)]_{\lambda=0}$.

Proof: We use local charts $U = \{z \in \mathbb{C}^n \mid |z_i| < 1\}$ for $Y$ and we may assume that $\sigma$ has support in $U$.

(1) In suitable coordinates we have

$$D(I)|U = (z_1 \ldots z_s)$$

$$D_0 \cup D(J)|U = (z_1 \ldots z_t) \quad (0 \leq t \leq s \leq n)$$

$$g(z) = A(z)\left|z_1^{v(1)} \ldots z_s^{v(s)}\right|^2 \quad (A(z) > 0, v(j) = v_{i(j)}).$$

We chose $l \in \mathbb{N}$ large enough to get smooth forms $(z_1 \ldots z_t)^l \psi_j$ on $U$. Then

$$M(\lambda) = \int_{\mathbb{C}^n} |z_1|^{2(v(1)\lambda-1)} \ldots |z_t|^{2(v(t)\lambda-1)} \left| z_{t+1}^{v(t+1)} \ldots z_s^{v(s)} \right|^{2\lambda} \sigma A^\lambda |z_1 \ldots z_t|^{2l} \psi_1 \wedge \overline{\psi_2}.$$

This integral is obtained from (2-1) by the substitution

$\lambda_i = v(i)\lambda - 1, 1 \le i \le t,$

$\lambda_i = v(i)\lambda, t+1 \le i \le s,$

$\lambda_i = 0, s+1 \le i \le n,$

$w = \lambda, \kappa(z,w)dV = \sigma A^w |z_1 \ldots z_t|^{2l}\psi_1 \wedge \overline{\psi_2}.$

Since the meromorphic continuation is obtained the same way, (1) follows.

(2) This reasoning shows, that the pole order of $M(\lambda)$ at 0 on the chart $U$ is at most $t \le k+1$.

(3) is proved similar to (1)

(4) By (1) and (2) it is enough to consider the case

$D_i|U = (z_i), 1 \le i \le k,$

$D_0|U = (z_{k+1})$

$D(I)|U = (z_1 \ldots z_s) \quad (k+1 \le s \le n)$

$g(z) = A(z)\left| z_1^{v(1)} \ldots z_s^{v(s)} \right|^2 \quad (A > 0, v(i) = v_i, 1 \le i \le k, v(k+1) = v_0)$

$\tilde{g}(z_{k+1},\ldots,z_n) = \tilde{A}(z_{k+1},\ldots,z_n)|z_{k+1}|^{2\tilde{v}} \quad (\tilde{A} > 0)$, and

$$\psi_j = z_{k+1}^{-l}\left[\frac{dz_1 \wedge \ldots \wedge dz_k}{z_1 \ldots z_k} \wedge \alpha_j + \beta_j\right], j = 1,2,$$

with forms $\alpha_j \in \mathcal{E}^{n-k}(U), \beta_j \in \sum_{L \subsetneq J} \frac{(dz)_L}{z_L} \wedge \mathcal{E}^{n-|L|}(U)$ and $l \in \mathbb{N}$.

The left-hand side results from

$$M(\lambda) = \int_{\mathbb{C}^n} |z_1|^{2(v_1\lambda-1)} \ldots |z_k|^{2(v_k\lambda-1)}|z_{k+1}|^{2(v_0\lambda-l)}\left| z_{k+2}^{v(k+2)} \ldots z_s^{v(s)}\right|^{2\lambda} \sigma A^\lambda (dz_1 \wedge \ldots \wedge dz_k \wedge$$
$$\alpha_1 + z_1 \ldots z_k\beta_1) \wedge \overline{(dz_1 \wedge \ldots \wedge dz_k \wedge \alpha_2 + z_1 \ldots z_k\beta_2)}.$$

If we substitute in (2-1)

$\lambda_i = v_i\lambda - 1, 1 \le i \le k,$

$\lambda_{k+1} = v_0\lambda - l,$

$\lambda_i = v(i)\lambda, k+2 \le i \le s,$

$\lambda_i = 0, s+1 \le i \le n,$

$w = \lambda, \kappa(z,w)dV = \sigma A^w \left|z_1 \ldots z_k z_{k+1}^l\right|^2 \psi_1 \wedge \overline{\psi_2},$

we obtain by Cor. 2.2:

$$[\lambda^{k+1}M(\lambda)]_{\lambda=0}v_0 \ldots v_k =$$
$$\pi^{k+1}\left(\frac{1}{(l-1)!}\right)^2 \int_{\{z|z_i=0,1\le i\le k+1\}}\left(L_{z_{k+1}}L_{\bar{z}_{k+1}}\right)^{l-1}\left[\sigma i\left(\frac{\partial}{\partial y_{k+1}}\right)i\left(\frac{\partial}{\partial x_{k+1}}\right)\ldots i\left(\frac{\partial}{\partial y_1}\right)i\left(\frac{\partial}{\partial x_1}\right)(dz_1 \wedge$$
$$\ldots \wedge dz_k \wedge \alpha_1 \wedge \overline{dz_1 \wedge \ldots \wedge dz_k \wedge \alpha_2})\right].$$

($x_i, y_i$ are the real coordinates, $L$ is the Lie-derivation and $i(w)$ the contraction with the vector field $w$.)

To evaluate the expression in brackets we need the constant

$$C(p) = (-2i)^p(-1)^{p(p-1)/2},$$

which satisfies $dz_1 \wedge \ldots \wedge dz_p \wedge \overline{dz_1 \wedge \ldots \wedge dz_p} = C(p)dV(z_1, \ldots, z_p)$. The restriction of this bracket to $\{z|z_i = 0, 1 \le i \le k\}$ reads

$$(-1)^{k(n-k)}C(k)\,\sigma i\left(\frac{\partial}{\partial y_{k+1}}\right)i\left(\frac{\partial}{\partial x_{k+1}}\right)\alpha_1 \wedge \overline{\alpha_2}\bigg|_{z_1=\cdots=z_k=0}.$$

Right-hand side: In our case $Res_{D_K}\psi_j = z_{k+1}^{-l}\alpha_j$, and

$$\widetilde{M}(\lambda) = \int_{\{z|z_i=0,1\le i\le k\}}|z_{k+1}|^{2(\tilde{v}\lambda-l)}\sigma\tilde{A}^\lambda\alpha_1 \wedge \overline{\alpha_2}.$$

By Cor. 2.2 we get

$$\left[\lambda\widetilde{M}(\lambda)\right]_{\lambda=0}\tilde{v} = \pi\left(\frac{1}{(l-1)!}\right)^2 \int_{\{z|z_i=0,1\le i\le k+1\}}\left(L_{z_{k+1}}L_{\bar{z}_{k+1}}\right)^{l-1}\left[\sigma i\left(\frac{\partial}{\partial y_{k+1}}\right)i\left(\frac{\partial}{\partial x_{k+1}}\right)(\alpha_1 \wedge \overline{\alpha_2})\right].$$

Putting both sides together, we obtain

$$[\lambda^{k+1}M(\lambda)]_{\lambda=0}v_0 \ldots v_k = \pi^k(-1)^{k(n-k)}C(k)\left[\lambda\widetilde{M}(\lambda)\right]_{\lambda=0}\tilde{v}.$$

### 2.2.2 Favorite situation II

In section 2.2.1 we compared the Mellin-type integral on a manifold to such an integral on a submanifold arising as intersection of certain divisors. Here we investigate further the latter one.

Let $N$ be a compact complex manifold of dimension $m$. Let $V \cup E$ be a reduced (effective) divisor with normal crossings and $V$ smooth. We write $N' := N\setminus E$, $V' := V\setminus E$. Let $h: N \to \mathbb{R}$ be a tube function along $V$ with multiplicity $v > 0$.

We consider closed semi-meromorphic forms

$$\xi_1 \in \Gamma(N, \mathcal{E}_N^m(*V)), \xi_2 \in \Gamma(N, \mathcal{E}_{N,E}^m(*V)),$$

and we want to determine $\left[\lambda\int_N h^\lambda \xi_1 \wedge \overline{\xi_2}\right]_{\lambda=0}$. Following [Ba 3], we reduce the pole order along $V$ to the logarithmic case.

**Lemma 2.5:** There are semi-meromorphic forms $\alpha \in \Gamma(N, \mathcal{E}_N^{m-1}(*V))$, $\beta \in \Gamma(N, \mathcal{E}_{N,E}^{m-1}(*V))$ such that $\xi_1' := \xi_1 + d\alpha \in \Gamma(N, \mathcal{E}_N^m(\log V))$, $\xi_2' := \xi_2 + d\beta \in \Gamma(N, \mathcal{E}_{N,E}^m(\log V))$.

Proof: The first assertion is [Ba 3, prop. 1]. The second one is obtained in a similar way. By Cor. A.1.6 the inclusion $\mathcal{E}_{N,E}^\bullet(\log V) \to \mathcal{E}_{N,E}^\bullet(*V)$ is a quasi-isomorphism of soft sheaves and induces an isomorphism $H^\bullet \Gamma \mathcal{E}^\bullet_{N,E}(\log V) \tilde{\to} H^\bullet \Gamma \mathcal{E}^\bullet_{N,E}(*V)$. The claim is, that this map is surjective in dimension $m$.

The Poincaré residue $Res_V \xi'_1 \in Z^{m-1} \Gamma \mathcal{E}^\bullet_V$ resp. $Res_{V'}(\xi'_1|N') \in Z^{m-1} \Gamma \mathcal{E}^\bullet_{V'}$ resp. $Res_V \xi'_2 \in Z^{m-1} \Gamma \mathcal{E}^\bullet_{V,E \cap V}$ defines a cohomology class in $H^{m-1}(V, \mathbb{C})$ resp. $H^{m-1}(V', \mathbb{C})$ resp. $H^{m-1}_c(V', \mathbb{C})$.

By Cor. A.1.6 (and the softness of the sheaves), these classes depend only on the class of $\xi_1$ resp. $\xi_1|N'$ resp. $\xi_2$ modulo exact semi-meromorphic forms. They are called Leray residues and will be denoted by

$$res_V \xi_1 = [Res_V \xi'_1], res_{V'} \xi_1 = [Res_{V'} \xi'_1|N'], res_{V',c} \xi_2 = [Res_V \xi'_2].$$

As is well known, $res_{V'} \xi_1$ is the image of $[\xi_1|N'\backslash V'] \in H^m(N'\backslash V', \mathbb{C})$ under the map $R$ in the Gysin sequence

(2-2) $\qquad H^m(N') \to H^m(N'\backslash V') \xrightarrow{R} H^{m-1}(V') \xrightarrow{\gamma} H^{m+1}(N')$

(i.e. the cohomology sequence of the lines in the diagram

$$\begin{array}{ccccccccc} 0 & \to & \Omega^\bullet_{N'} & \to & \Omega^\bullet_{N'}(\log V') & \to & \Omega^{\bullet-1}_{V'} & \to & 0 \\ & & \downarrow \approx & & \downarrow \approx & & \downarrow \approx & & \\ 0 & \to & K^\bullet & \to & \mathcal{E}^\bullet_{N'}(\log V') & \to & \mathcal{E}^{\bullet-1}_{V'} & \to & 0 \end{array}$$

where $K^\bullet \supseteq \mathcal{E}^\bullet_{N'}$ is the kernel, cf. A.1.6 and [Gr-S]). If we take the Gysin sequence for $V \subseteq N$, we have the analogous representation for $res_V \xi_1$, and $res_{V'} \xi_1$ is the image of $res_V \xi_1$ under the restriction $H^{m-1}(V) \to H^{m-1}(V')$.

Our aim is now to give a similar representation for $res_{V',c} \xi_2$.

We need some more notation. We have the inclusions:

$\qquad j: N'\backslash V' \to N$

$\qquad j_0: N\backslash V \to N$

$\qquad j_1: N' \to N$

$\qquad j': N'\backslash V' \to N\backslash V$

$\qquad i: V' \to N'$

We define supporting families for the spaces $N, N', N'\backslash V', V, V'$:

$\qquad c_N := \{A \subseteq N | A \text{ compact}\}$

$\qquad c_{N'} := c_N | N' = \{A \in c_N | A \subseteq N'\}$

$\qquad \phi := c_N \cap (N\backslash V) | N'\backslash V' = \{B := A \cap (N\backslash V) | A \in c_N, B \subseteq N\backslash V\backslash E\}$
$\qquad \qquad \qquad = \{B \subseteq N\backslash V\backslash E | B \text{ closed in } N\backslash V\}$

$\qquad c_V := c_N \cap V = \{B \subseteq V | B \text{ compact}\}$

$\qquad c_{V'} := c_V | V'$

They are paracompactifying in the sense of [Go, p. 150].[5]

**Lemma 2.6:**

(1) $H_c^k(N') := H_{c_{N'}}^k(N', \mathbb{C}) = \mathbb{H}^k \mathcal{E}_{N,E}^\bullet = \mathbb{H}^k \Omega_{N,E}^\bullet$ and
$H_c^k(V') := H_{c_{V'}}^k(V', \mathbb{C}) = \mathbb{H}^k \mathcal{E}_{V,E\cap V}^\bullet = \mathbb{H}^k \Omega_{V,E\cap V}^\bullet$.

(2) $H_\phi^k(N'\backslash V', \mathbb{C}) = \mathbb{H}^k j_{0*} j_0^{-1} \mathcal{E}_{N,E}^\bullet = \mathbb{H}^k \mathcal{E}_{N,E}^\bullet(*V) = \mathbb{H}^k \mathcal{E}_{N,E}^\bullet(\log V) = \mathbb{H}^k \Omega_{N,E}^\bullet(\log V) = \mathbb{H}^k \Omega_{N,E}^\bullet(*V)$.

Proof: (1) By [Go, p. 190] (and since $c_N = \{A \subseteq N | A \text{ closed}\}$) we have

$$H_{c_N|N'}^k(N', \mathbb{C}) = H_{c_N}^k(N, j_{1!}\mathbb{C}) = H^k(N, j_{1!}\mathbb{C}).$$

Because of $j_{1!}\mathbb{C} \approx \Omega_{N,E}^\bullet \approx \mathcal{E}_{N,E}^\bullet$ (Lemma A.1.1), this is isomorphic to $\mathbb{H}^k \mathcal{E}_{N,E}^\bullet$ and $\mathbb{H}^k \Omega_{N,E}^\bullet$. The second assertion is of the same kind.

(2) Again, we have ([Go, p. 190])

$$H_\phi^k(N\backslash V\backslash E, \mathbb{C}) = H_{c_N \cap (N\backslash V)}^k(N\backslash V, j_!'\mathbb{C}) = H^k(N\backslash V, j_!'\mathbb{C}),$$

and this is isomorphic to

$$H^k \Gamma(N\backslash V, \mathcal{E}_{N,E}^\bullet) = \mathbb{H}^k j_{0*} j_0^{-1} \mathcal{E}_{N,E}^\bullet,$$

where we use, that $j_{0*} j_0^{-1} \mathcal{E}_{N,E}^\bullet$ is soft (and fine), being a module over a soft (and fine) sheaf of rings ([Go, p. 156-157]). The rest of (2) follows, because the complexes

$$j_{0*} j_0^{-1} \mathcal{E}_{N,E}^\bullet, \mathcal{E}_{N,E}^\bullet(*V), \mathcal{E}_{N,E}^\bullet(\log V), \Omega_{N,E}^\bullet(\log V), \Omega_{N,E}^\bullet(*V), j_{0*} j_0^{-1} \Omega_{N,E}^\bullet$$

are all quasi-isomorphic (A.1.4, A.1.6, A.1.10).

From the lines in the diagram

(2-3) $$0 \to \Omega_{N,E}^\bullet \to \Omega_{N,E}^\bullet(\log V) \xrightarrow{Res} \Omega_{V,E\cap V}^{\bullet-1} \to 0$$
$$\downarrow\approx \quad\quad \downarrow\approx \quad\quad \downarrow\approx$$
$$0 \to K^\bullet \to \mathcal{E}_{N,E}^\bullet(\log V) \xrightarrow{Res} \mathcal{E}_{V,E\cap V}^{\bullet-1} \to 0$$

($K^\bullet \supseteq \mathcal{E}_{N,E}^\bullet$ the kernel) arises the sequence

(2-4) $$\to H_c^m(N') \to H_\phi^m(N\backslash V\backslash E, \mathbb{C}) \xrightarrow{R_c} H_c^{m-1}(V') \xrightarrow{\gamma_c} H_c^{m+1}(N').$$

We now have:

 $res_{V',c} \xi_2$ *is the image of* $[\xi_2]$ *under* $R_c$.

---

[5] Proof for $\phi$: $N\backslash V$ is paracompact, hence normal. Therefore, the closed and disjoint subsets $E\backslash V$ and $B$ of $N\backslash V$ have disjoint neighborhoods $U_1, U_2$ in $N\backslash V$. Then $\overline{U_2} \in \phi$.

Furthermore is clear, that $res_{V',c}\xi_2$ is mapped to $res_V\xi_2$ and $res_{V'}\xi_2$ by the canonical maps $H_c^{m-1}(V') \to H^{m-1}(V) \to H^{m-1}(V')$.

We shall give still another more suggestive representation for (2-4).

**Remark 2.7:** The canonical map

$$H^k\Gamma_{c_{N'}}(N', \mathcal{E}_{N'}^\bullet(\log V')) \to H^k\Gamma(N, \mathcal{E}_{N,E}^\bullet(\log V))$$

is an isomorphism.

Proof: The left-hand side can be rewritten as

$$H^k\Gamma_{c_{N'}}(N', \mathcal{E}_{N'}^\bullet(\log V')) = H^k\Gamma(N, j_{1!}\mathcal{E}_{N'}^\bullet(\log V')) = \mathbb{H}^k(N, j_{1!}\mathcal{E}_{N'}^\bullet(\log V'))$$
$$= \mathbb{H}^k(N, j_{1!}\Omega_{N'}^\bullet(\log V')) = \mathbb{H}^k(N, j_{1!}j_1^{-1}\Omega_{N,E}^\bullet(\log V))$$

(since $j_{1!}$ is exact and transforms soft sheaves into such [Go, p. 154]).

Because of the exact sequence [Go, p. 140]

$$0 \to j_{1!}j_1^{-1}\Omega_{N,E}^\bullet(\log V) \to \Omega_{N,E}^\bullet(\log V) \to \Omega_{N,E}^\bullet(\log V)|E \to 0,$$

it is enough to show that $\forall\, k \geq 0$

$$\mathbb{H}^k(E, \Omega_{N,E}^\bullet(\log V)|E) = 0.$$

Because of

$$0 \to \Omega_{N,E}^\bullet|E \to \Omega_{N,E}^\bullet(\log V)|E \to \Omega_{V,E\cap V}^{\bullet-1}|E \to 0,$$

it is sufficient if $\forall\, k \geq 0$

$$\mathbb{H}^k(E, \Omega_{N,E}^\bullet|E) = 0,\ \mathbb{H}^k(E, \Omega_{V,E\cap V}^{\bullet-1}|E) = 0.$$

This is true since the complexes are resolutions of zero (A.1.1).

We see from this remark, that (2-4) is obtained also like this: One applies $\Gamma_{c_{N'}}(N',*)$ to

$$0 \to K^\bullet \to \mathcal{E}_{N'}^\bullet(\log V') \to \mathcal{E}_{V'}^{\bullet-1} \to 0$$

(which remains exact by softness of the sheaves) and forms the exact cohomology sequence.

**Lemma 2.8:** The maps $i^*: H^{m-1}(N') \to H^{m-1}(V')$ and $\gamma_c/(-2\pi i): H_c^{m-1}(V') \to H_c^{m+1}(N')$ are adjoint with respect to Poincaré duality ($H^{m-1}(V') \cong H_c^{m-1}(V')^*$, $H^{m-1}(N') \cong H_c^{m+1}(N')^*$).

Proof: This is analogous to [Gr-S, 2.18]. We have to show $\langle \gamma_c x, y\rangle = -2\pi i \langle x, i^*y\rangle$ for $x \in H_c^{m-1}(V')$, $y \in H^{m-1}(N')$. We determine $\gamma_c x$ from the lower line in (2-3). Without loss of generality, let $x = [\varphi]$ with $\varphi \in \Gamma_c(\mathcal{E}_{V'}^{m-1})$ and $y = [\psi] \in \Gamma(\mathcal{E}_{N'}^{m-1})$. Consider a tube function $g$ for $V \subseteq N$, i.e. locally

$$g = A|t|^2,\ \frac{\partial g}{g} = \partial \log g = \frac{dt}{t} + \frac{\partial A}{A}.$$

By means of a retraction to $V$, $\varphi$ can be extended to $\tilde{\varphi} \in \Gamma_c(\mathcal{E}_{N'}^{m-1})$ such that $d\tilde{\varphi} = 0$ in some neighborhood of $V$ in $N$. Then $\phi := d\left(\frac{\partial g}{g} \wedge \tilde{\varphi}\right) = \bar{\partial}\partial \log g \wedge \tilde{\varphi} - \frac{\partial g}{g} \wedge d\tilde{\varphi} \in \Gamma_c(\mathcal{E}_{N'}^{m+1})$ is a representative of $\gamma_c x$, and by a simple local calculation (using Cauchy's formula)

$$\langle \gamma_c x, y\rangle = \int_N \phi \wedge \psi = \lim_{\varepsilon \to 0}\int_{g\geq \varepsilon} d\left(\frac{\partial g}{g} \wedge \tilde{\varphi} \wedge \psi\right) = -\lim_{\varepsilon \to 0}\int_{g=\varepsilon} \frac{\partial g}{g} \wedge \tilde{\varphi} \wedge \psi = -2\pi i \int_V \tilde{\varphi} \wedge \psi = -2\pi i \langle x, i^* y\rangle.$$

Next, we compare $\int_N h^\lambda \xi_1 \wedge \bar{\xi}_2$ to $\int_N h^\lambda \xi_1' \wedge \bar{\xi}_2'$.

**Lemma 2.9:** As functions of $\lambda$ (with the notation of Lemma 2.5)

$$\int_N h^\lambda \xi_1 \wedge \bar{\xi}_2, \int_N h^\lambda \frac{dh}{h} \wedge \alpha \wedge \bar{\xi}_2, \int_N h^\lambda \frac{dh}{h} \wedge \xi_1 \wedge \bar{\beta}, \int_N h^\lambda \frac{dh}{h} \wedge \alpha \wedge d\bar{\beta}$$

are meromorphic (extendable) with at most simple poles ($\in \frac{1}{v}\mathbb{Z}$).

Proof: By a partition of unity, we obtain summands of the form

$$\int_U (A|z_1|^{2v})^\lambda |z_1|^{-2l} \sigma(z) dV, \sigma \in C_c^\infty(U), l \geq 0,$$

on the chart $U \subseteq \mathbb{C}^m$ with $V|U = (z_1)$, $h|U = A|z_1|^{2v}$. We apply Cor. 2.2 (with $\kappa(z,w) = A^w \sigma(z)$).

**Lemma 2.10:** We have $\left[\lambda \int_N h^\lambda \xi_1 \wedge \bar{\xi}_2\right]_{\lambda=0} = \left[\lambda \int_N h^\lambda \xi_1' \wedge \bar{\xi}_2'\right]_{\lambda=0}$.

Proof: We use Stokes' theorem in the form

$$0 = \int_N d(h^\lambda \alpha \wedge \bar{\xi}_2) = \lambda \int_N h^\lambda \frac{dh}{h} \wedge \alpha \wedge \bar{\xi}_2 + \int_N h^\lambda d\alpha \wedge \bar{\xi}_2,$$

$$0 = \int_N d(h^\lambda \xi_1 \wedge \bar{\beta}) = \lambda \int_N h^\lambda \frac{dh}{h} \wedge \xi_1 \wedge \bar{\beta} \pm \int_N h^\lambda \xi_1 \wedge d\bar{\beta},$$

$$0 = \int_N d(h^\lambda \alpha \wedge d\bar{\beta}) = \lambda \int_N h^\lambda \frac{dh}{h} \wedge \alpha \wedge d\bar{\beta} + \int_N h^\lambda d\alpha \wedge d\bar{\beta},$$

where $\operatorname{Re} \lambda$ has to be chosen large enough (e.g. $\operatorname{Re} \lambda > 1 + l$), such that on the left occur the derivatives of $C^1$ forms. The integrals in 2.10 differ by the second summands on the right of these equations. By 2.9 the first summands do not have a pole at $\lambda = 0$.

**Lemma 2.11:** $\left[\lambda \int_N h^\lambda \xi_1' \wedge \bar{\xi}_2'\right]_{\lambda=0} v = (2\pi i)(-1)^m \int_V \operatorname{Res}_V \xi_1' \wedge \overline{\operatorname{Res}_V \xi_2'}$.

Proof: Let $U \subseteq \mathbb{C}^m$ be a local chart for $N$ with $V|U = (z_1)$, $h|U = A|z_1|^{2v}$, $\xi_j' = \frac{dz_1}{z_1} \wedge \alpha_j + \beta_j$, $\alpha_j, \beta_j \in \Gamma(U, \mathcal{E}_U^\bullet)$, $j = 1, 2$. For $\sigma \in C_c^\infty(U)$ the contribution to the left integral is

$$\int \sigma h^\lambda \xi_1' \wedge \bar{\xi}_2' = \int |z_1|^{2(v\lambda - 1)} \sigma A^\lambda (dz_1 \wedge \alpha_1 + z_1 \beta_1) \wedge \overline{(dz_1 \wedge \alpha_2 + z_1 \beta_2)}.$$

If we substitute in (2-1)

$$\lambda_1 = v\lambda - 1, \lambda_i = 0, 2 \leq i \leq n = m,$$

$$w = \lambda \text{ and } \kappa(z,w)dV = \sigma A^w |z_1|^2 \xi_1' \wedge \overline{\xi_2'},$$

we obtain by cor. 2.2

$$\left[\lambda \int \sigma h^\lambda \xi_1' \wedge \overline{\xi_2'}\right]_{\lambda=0} v = \pi \int_{z_1=0} \sigma i\left(\frac{\partial}{\partial y_1}\right) i\left(\frac{\partial}{\partial x_1}\right) \left[(dz_1 \wedge \alpha_1 + z_1\beta_1) \wedge \overline{(dz_1 \wedge \alpha_2 + z_1\beta_2)}\right] =$$
$$\pi(-1)^{m-1}(-2i) \int_{z_1=0} \sigma \alpha_1 \wedge \overline{\alpha_2}.$$

**Lemma 2.12:** $\left[\lambda \int_N h^\lambda \xi_1 \wedge \overline{\xi_2}\right]_{\lambda=0} v = (2\pi i)(-1)^m \int_{V'} res_{V'} \xi_1 \cup \overline{res_{V',c}\xi_2},$

where $\cup: H^{m-1}(V') \times H_c^{m-1}(V') \to H_c^{2m-2}(V')$ ist the cup-product.

Proof: With $r: H^{m-1}(V) \to H^{m-1}(V')$ and $c: H_c^{m-1}(V') \to H^{m-1}(V)$ meaning the canonical maps, we have $r(res_V \xi_1) = res_{V'} \xi_1$, $c(res_{V',c}\xi_2) = res_V \xi_2$. If we write $\cup$ also for the cup-product on $H^{m-1}(V)$, we have the compatibility

$$\int_{V'} res_{V'} \xi_1 \cup \overline{res_{V',c}\xi_2} = \int_V res_V \xi_1 \cup \overline{res_V \xi_2}.$$

### 2.2.3 Conclusion

We are going to apply Lemma 2.12 to the situation 2.2.1 by setting: $m = n - k$,

$$N = D_J, V := Z \cap D_J, E := D' \cap D_J \text{ with } D' := D_{k+1} \cup \ldots \cup D_K$$

$$\xi_j := Res_{D_J} \psi_j, j = 1,2, \sigma \in C_c^\infty(Y), \sigma = 1 \text{ near } Z \cap D_J.$$

Here we additionally suppose: $D_J$ is compact, and $\psi_j \in \Gamma(Y, \mathcal{E}_{Y,D}^n,(\log D(J))(*Z))$ is closed, hence also $\xi_j \in \Gamma(N, \mathcal{E}_{N,E}^m(*V)), j = 1,2$. In particular, these forms define both cohomology classes of compact support.

As already observed, we need not specify the order in which the iterated Poincaré residue $Res_{D_J}\psi_j$ along the divisors $D_i, i \in J$, is taken, because this cancels out.

By Lemma 2.4, (4), it follows for the function $M(\lambda) = \int_Y g^\lambda \sigma \psi_1 \wedge \overline{\psi_1}$ of section 2.2.1:

**Corollary 2.13:**

$$P_{-(k+1)}(0, M) = (v_0 \ldots v_k)^{-1}(2\pi i)^{k+1}(-1)^{(k+1)(2n-k)/2} \int_{V'} (res_{V'} \xi_1) \cup \overline{(res_{V',c}\xi_2)}.$$

Later, we shall make essential use of the following remark.

**Remark 2.14**: By our assumptions, $res_{V',c}\xi_1$ is defined, and $res_{V'}\xi_1$ is the image under $c: H_c^{m-1}(V') \to H^{m-1}(V')$. Therefore $\int_{V'} (res_{V'} \xi_1) \cup \overline{(res_{V',c}\xi_2)}$ depends only on the image $res_{V'}\xi_2$ of $res_{V',c}\xi_2$, and we write for it

$$\langle res_{V'}\xi_1, \overline{res_{V'}\xi_2}\rangle_{V'}.$$

Proof: There is a well-defined bilinear form on $image(c) \times image(c)$.

## 3   Application to $\int |f|^{2\lambda} \square$ and the coefficients of asymptotic expansions

### 3.1   Integer degrees

**3.1.1** We go back to the setting of section 0.4 with an integer $a \in \mathbb{N}^*$ and $1 \leq r \leq n-1$. We consider the integral

$$M(\lambda) = \int_X |f|^{2\lambda} \sigma \varphi_1 \wedge \overline{\varphi_2}, \quad (\sigma \in C_c^\infty(X), \sigma = 1 \text{ near } 0)$$

with $\varphi_1, \varphi_2$ as in section 0.4. We aim to show that $(\alpha, k)(\varphi_1) = (a-1, r)$, and for this it is enough, if we can choose $\varphi_2$ such that

(3-1)    $P_{-(r+1)}(-a, M) \neq 0$.

We continue to denote by $X_\Sigma{'}$ the inverse image of $X$ under $\rho_2 \rho_1 : X_\Sigma \to \mathbb{C}^n$. As in section 1, $\bigcup_{\check\tau \in S} X_\tau$ is the open toric subset, where $f^{-a} \varphi_1$ has $r+1$ "pole components" ($Z$ and $r$ exceptional divisors). $\overline{T}_\tau \subseteq X_\Sigma$ is the closure of the stratum $T_\tau \subseteq X_\tau$, and $\overline{T}_\tau \subseteq U_\tau \subseteq X_\Sigma{'}$ are disjoint neighborhoods.

Outside of $\bigcup_{\check\tau \in S} U_\tau$, the forms $f^{-a} \varphi_j$ have less than $r+1$ "pole components", i.e. they are in $W_{r-1} \Omega^n(\log D)(aZ)$.[6] As a consequence

$$P_{-(r+1)}(-a, M) = \sum_{\check\tau \in S} P_{-(r+1)}\left(\lambda = 0, \int_{U_\tau} |f|^{2\lambda} \sigma_\tau \frac{\varphi_1}{f^a} \wedge \frac{\overline{\varphi_2}}{\overline{f^a}}\right),$$

where $\sigma_\tau \in C_c^\infty(U_\tau)$ and $= 1$ near $\overline{T}_\tau$ (i.e. the leading Laurent coefficient depends only on $\bigcup_{\check\tau \in S} U_\tau$, cf. section 2.1).

We apply Cor. 2.13 in section 2.2.3 to the right-hand side. For $\check\tau \in S$ we put

$Y := U_\tau$

$D(J) := \bigcup_{l \leq \check\tau} D(l) \ (\Rightarrow k = r), D_0 := Z,$

$D' := \bigcup_{\substack{l \in \Sigma^{(1)} \\ l \not\leq \check\tau}} D(l) \ (\Rightarrow D(J) \cup D' = X_\Sigma \setminus (\mathbb{C}^*)^n),$

$g := |f|^2, \sigma := \sigma_\tau \ (\in C_c^\infty(U_\tau), = 1 \text{ near } \overline{T}_\tau),$

$\psi_j = \varphi_j / f^a, j = 1,2.$

Then we get (cf. Remark 2.14)

---

[6] In the notation of section 1: $f^{-a} \varphi_j$ has only poles along $Z$ and $D(l)$, $l \in \Sigma^{(1)}$ (1-skeleton) with $l \subseteq \check\sigma(\delta)$. The set of all non-empty intersections of $r$ such $D(l)$ is just $\overline{T}_\tau = \bigcap_{l \leq \check\tau} D(l), \check\tau \subseteq \check\sigma(\delta)$ cone of dimension $r$.

$$P_{-(r+1)}\left(\lambda = 0, \int_{U_\tau} |f|^{2\lambda} \sigma_\tau \frac{\varphi_1}{f} \wedge \overline{\frac{\varphi_2}{f}}\right) = (2\pi i)^{r+1}\varepsilon(n,r)v_\tau^{-1} \langle res\left(\frac{\varphi_1}{f^a}\right), \overline{res\left(\frac{\varphi_2}{f^a}\right)}\rangle_{T_\tau \cap Z},$$

with

$$res\left(\frac{\varphi_j}{f^a}\right) := res_{T_\tau \cap Z} Res_{T_\tau}\left(\frac{\varphi_j}{f^a}\right) \in H^{m-1}(T_\tau \cap Z, \mathbb{C}) \ (m = n - r),$$

$\varepsilon(n,r) = \pm 1$ (independent of $\check{\tau} \in S$) and $v_\tau \in \mathbb{N}^+$ (depending on the multiplicities of $f$ along $D$).

By Lemma 1.3, these terms differ for $\check{\tau} \in S$ only by positive factors (since the non-zero factor in 1.3 occurs twice). Therefore it is sufficient to fix some $\check{\tau} \in S$. We have to find a form $\varphi_2 \in \Omega_0^n$, $supp(\varphi_2) \subseteq a\delta°$, such that

(3-2) $\quad \langle res\left(\frac{\varphi_1}{f^a}\right), \overline{res\left(\frac{\varphi_2}{f^a}\right)}\rangle_{T_\tau \cap Z} \neq 0.$

Write $F_{\check{\tau}} := \overline{T}_\tau$, $D_\tau := \overline{T}_\tau \setminus T_\tau$, $Z_\tau := Z \cap F_{\check{\tau}}$, $V_\tau := Z \cap T_\tau$, $i: V_\tau \hookrightarrow T_\tau$, and

$$res_c\left(\frac{\varphi_j}{f^a}\right) := res_{V_\tau,c} Res_{F_{\check{\tau}}}\left(\frac{\varphi_j}{f^a}\right) \in H_c^{m-1}(V_\tau, \mathbb{C})$$

as in section 2.2.2. Also, recall the exact sequence

$$H_c^m(T_\tau) \to H_\phi^m(T_\tau \setminus V_\tau) \xrightarrow{R_c} H_c^{m-1}(V_\tau) \xrightarrow{\gamma_c} H_c^{m+1}(T_\tau).$$

**Lemma 3.1:**

(1) $\left\{Res_{F_{\check{\tau}}}\left(\frac{\varphi}{f^a}\right) \middle| supp(\varphi) \subseteq a\delta°\right\} = \Gamma\left(F_{\check{\tau}}, \Omega_{F_{\check{\tau}}}^m(aZ_\tau)\right)$

(2) For $a > m$ the map $res_c: \Gamma\left(F_{\check{\tau}}, \Omega_{F_{\check{\tau}}}^m(aZ_\tau)\right) \to \ker \gamma_c$ is surjective.

Proof: Follows from A.2.2 and A.2.6.

According to Lemma 2.8, $i^*: H^{m-1}(T_\tau) \to H^{m-1}(V_\tau)$ and $\gamma_c/(-2\pi i)$ are adjoint with respect to Poincaré duality. This implies:

**Lemma 3.2:** For $x \in H^{m-1}(V_\tau)$ we have: $\langle x, y \rangle = 0 \ \forall \ y \in \ker \gamma_c \Leftrightarrow x \in image(i^*)$.

Proof: $\langle x, y \rangle = 0 \ \forall \ y \in \ker \gamma_c \Leftrightarrow \langle x, . \rangle$ defines a linear form on $image(\gamma_c)$, which can be written as $\langle z, . \rangle$ for some $z \in H^{m-1}(T_\tau) \Leftrightarrow \langle x, y \rangle = \langle z, \gamma_c(y) \rangle = (-2\pi i)\langle i^*(z), y \rangle \ \forall y \in H_c^{m-1}(V_\tau) \Leftrightarrow x = (-2\pi i)i^*(z)$.

Let $H_c^{m-1}(V_\tau) \xrightarrow{c} H^{m-1}(V_\tau) \xleftarrow{i^*} H^{m-1}(T_\tau)$ be the canonical maps. By A.2.8 we have:

**Lemma 3.3:** For $m \geq 2$, $image(c) \cap image(i^*) = 0$ in $H^{m-1}(V_\tau)$, for $m = 1$, $image(R_c) \cap image(i^*) = 0$ in $H^0(V_\tau)$.

**Lemma 3.4:** The intersection of the images of the canonical maps

$$H^m(T_\tau) \to H^m(T_\tau \setminus V_\tau) \leftarrow H^m_\phi(T_\tau \setminus V_\tau)$$

is zero (for $m \geq 1$).

Proof: See A.2.9.

We apply these lemmas to $x = res(\frac{\varphi_1}{f^a}) \in image(c) \subseteq H^{m-1}(V_\tau)$. If $\left[Res_{T_\tau}(\frac{\varphi_1}{f^a})\right] \neq 0$ in $H^m(T_\tau \setminus V_\tau)$, then $x \neq 0$ (by 3.4 and (2-2)), hence $x \notin image(i^*)$ by 3.3, and there is a $y \in \overline{\ker \gamma_c}$ with $\langle x, y \rangle \neq 0$ by 3.2. In the case $a > m$ this is representable as $y = \overline{res_c(\frac{\varphi_2}{f^a})}$ with $supp(\varphi_2) \subseteq a\delta°$ (by 3.1, since the conjugate $\bar{y} \in \ker \gamma_c$).

We thus have proved for $a > m$:

**Proposition 3.5:** Suppose $a \in \mathbb{N}^+$, $1 \leq r \leq n-1$, and $\left[Res_{T_\tau}(\frac{\varphi_1}{f^a})\right] \neq 0$ in $H^{n-r}(T_\tau \setminus V_\tau)$ (for one and thereby all $\check{\tau} \in S$). Then $(\alpha, k)(\varphi_1) = (a-1, r-1)$.

Proof: There remains the case $a \leq m$. But we have for $b \geq 0$

$$Res_{T_\tau}\left(\frac{\varphi_1}{f^a}\right) = Res_{T_\tau}\left(\frac{f^b \varphi_1}{f^{a+b}}\right) = Res_{T_\tau}\left(\frac{f_\delta^b \varphi_1}{f^{a+b}}\right)$$

(since $(f^b - f_\delta^b)\varphi_1$ has Newton order $> a+b$) and $supp(f_\delta^b \varphi_1) \subseteq (a+b)\delta°$. If we choose $a + b > n - r$ we are in the case already proved and get $(\alpha, k)(f_\delta^b \varphi_1) = (a + b - 1, r - 1)$. Then $(\alpha, k)(f^b \varphi_1) = (a + b - 1, r - 1)$, because of $(\alpha, k)\left((f^b - f_\delta^b)\varphi_1\right) \geq (v, l)\left((f^b - f_\delta^b)\varphi_1\right) - (1, 0) > (a + b - 1, r - 1)$. Finally $(\alpha, k)(\varphi_1) = (\alpha, k)(f^b \varphi_1) - (b, 0) = (a - 1, r - 1)$.

### 3.1.2 Evaluation of the condition $0 \neq \left[Res_{T_\tau}(\frac{\varphi_1}{f^a})\right] \in H^{n-r}(T_\tau \setminus V_\tau)$

Consider again $\check{\tau} \in S$, $T_\tau \subseteq X_\tau$ the corresponding stratum of dimension $n - r$, and $V_\tau = Z \cap T_\tau$, $F_{\check{\tau}} = \bar{T}_\tau \subseteq X_\Sigma$, $Z_\tau = Z \cap F_{\check{\tau}}$.

Recall the representation of $Res_{T_\tau}(\frac{\varphi_1}{f^a})$ from Lemma 1.3. With $m \in \delta \cap \mathbb{Z}^n$, $\tilde{f} = x^{-m}f$, $\tilde{\varphi}_1 := x^{-am}\varphi_1$, this is the image

$$f^{-a}\varphi_1 = \tilde{f}^{-a}\tilde{\varphi}_1 \in \tilde{f}^{-a}\mathbb{C}[\tau \cap \mathbb{Z}^n] \otimes \Lambda^n \mathbb{C}^n$$
$$\downarrow \qquad\qquad \downarrow r \otimes I$$
$$\xi \in \tilde{f}|_{T_\tau}^{-a}\mathbb{C}[cospan(\tau) \cap \mathbb{Z}^n] \otimes \Lambda^{n-r}\mathbb{C}^n$$

under the map, which is composed of the restriction $r: \mathbb{C}[\tau \cap \mathbb{Z}^n] \to \mathbb{C}[cospan(\tau) \cap \mathbb{Z}^n]$ and the iteration $I = \Lambda_l i(l): \Lambda^n \mathbb{C}^n \to \Lambda^{n-r}\mathbb{C}^n$ of the contractions $i(l)$ ($l$ the primitive edge vectors of $\check{\tau}$). We note that $I$ maps to $\Lambda^{n-r} \bigcap_l \ker(l)$.

In order to describe the lower dimensional torus embedding $F_{\check{\tau}}$, we introduce new notation (valid in this context):

$$M := \mathbb{Z}^n$$

$$M' := M \cap \bigcap_{j \in J} \ker(l_j) \; (l_j, j \in J, \text{ the primitive edge vectors of } \check{\tau})$$

$$T' := \operatorname{Spec} \mathbb{C}[M'] = T_\tau, \; X' := F_{\check{\tau}} \supseteq T', \; D' := X' \backslash T'$$

$$g := \tilde{f}_{|T_\tau} \in \mathbb{C}[M'], \; Z' := \overline{(g|T')} \subseteq X' \; (\Rightarrow x^m g = f_\delta, \; Z' = Z_\tau)$$

$$\Delta' := \Delta(g) \subseteq M'_\mathbb{R} \text{ the Newton polyhedron}$$

Then $\xi \in \Gamma(X', \Omega^{n-r}(aZ')) \subseteq \Gamma(X', \Omega^{n-r}(\log D')(aZ'))$ and (Lemma A.2.6)

$$[\xi] \in \frac{\Gamma(X', \Omega^{n-r}(\log D')(aZ'))}{d\Gamma(X', \Omega^{n-r-1}(\log D')((a-1)Z'))} \subseteq H^{n-r}(T' \backslash Z').$$

The two vector spaces in the formula have representations as (Lemma A.2.2)

$$\Gamma(X', \Omega^{n-r}(\log D')(aZ')) = g^{-a} L(a\Delta') \otimes \Lambda^{n-r} M'_\mathbb{C},$$

$$\Gamma(X', \Omega^{n-r-1}(\log D')((a-1)Z')) = g^{-a+1} L((a-1)\Delta') \otimes \Lambda^{n-r-1} M'_\mathbb{C}.$$

The condition $[\xi] \neq 0$ reads after multiplication by $g^a$:

$$[g^a \xi] \neq 0 \quad \text{in} \quad \frac{L(a\Delta') \otimes \Lambda^{n-r} M'_\mathbb{C}}{\{-(a-1) dg \wedge \alpha + g d\alpha | \alpha \in L((a-1)\Delta') \otimes \Lambda^{n-r-1} M'_\mathbb{C}\}}$$

By the iterated contraction $id \otimes \Lambda_j i(l_j) : \mathbb{C}[M'] \otimes \Lambda^n M_\mathbb{C} \to \mathbb{C}[M'] \otimes \Lambda^{n-r} M'_\mathbb{C}$,

$$\{-(a-1) dg \wedge \beta + g d\beta | \beta \in L((a-1)\Delta') \otimes \Lambda^{n-1} M_\mathbb{C}\}$$

is mapped onto

$$\{-(a-1) dg \wedge \alpha + g d\alpha | \alpha \in L((a-1)\Delta') \otimes \Lambda^{n-r-1} M'_\mathbb{C}\}.$$

Proof: The first set is mapped into the second. This follows from $i(l_j) dg = 0$ and $L(l_j)\beta = (di(l_j) + i(l_j)d)\beta = 0$ for $\beta \in L((a-1)\Delta') \otimes \Lambda^p M_\mathbb{C}$ (observing $M' \subseteq \ker(l_j)$), so that $i(l_j)(dg \wedge \beta) = -dg \wedge i(l_j)\beta$ and $i(l_j) d\beta = L(l_j)\beta - di(l_j)\beta = -di(l_j)\beta$. By choosing $\beta = (\Lambda_j v_j) \wedge \alpha$, where $v_j, j \in J$, is a basis for a complement of $M'_\mathbb{C}$ in $M_\mathbb{C}$, we see also that the map is surjective.

So far, we have obtained the equivalent condition: $[\tilde{\varphi}_1] \neq 0$ in

$$L(a\Delta') \otimes \Lambda^n M_\mathbb{C} / \{-(a-1) dg \wedge \beta + g d\beta | \beta \in L((a-1)\Delta') \otimes \Lambda^{n-1} M_\mathbb{C}\}.$$

We still have to convert the denominator:

$$x^{am}\{-(a-1)dg \wedge \beta + gd\beta | \beta \in L((a-1)\Delta') \otimes \Lambda^{n-1}M_{\mathbb{C}}\} = df_\delta \wedge d(x^{(a-1)m}L((a-1)\Delta') \otimes \Lambda^{n-2}M_{\mathbb{C}}.$$

Proof: By the identity $x^{am}(-(a-1)dg \wedge \beta + gd\beta) = -(a-1)d(x^m g) \wedge x^{(a-1)m}\beta + x^m g d(x^{(a-1)m}\beta)$, the left-hand side is

$$\{-(a-1)df_\delta \wedge \alpha + f_\delta d\alpha | \alpha \in x^{(a-1)m}L((a-1)\Delta') \otimes \Lambda^{n-1}M_{\mathbb{C}}\}.$$

Inclusion $\supseteq$: For $df_\delta \wedge d\gamma$, $\gamma \in x^{(a-1)m}L((a-1)\Delta') \otimes \Lambda^{n-2}M_{\mathbb{C}}$ we can choose $\alpha := d\gamma \in x^{(a-1)m}L((a-1)\Delta') \otimes \Lambda^{n-1}M_{\mathbb{C}}$ with $d\alpha = 0$.

Inclusion $\subseteq$: By quasi-homogeneity, there is a vector field $v \in \tilde{M}_{\mathbb{C}}$ with $v(f_\delta) = f_\delta$. For $\alpha \in x^{(a-1)m}L((a-1)\Delta') \otimes \Lambda^{n-1}M_{\mathbb{C}}$ we then have the identity

$$-(a-1)df_\delta \wedge \alpha + f_\delta d\alpha = df_\delta \wedge (-(a-1)\alpha + i_v d\alpha) = -df_\delta \wedge di_v \alpha.$$

To see this, we use $df_\delta \wedge d\alpha = 0$ and $L_v \alpha = (a-1)\alpha$, obtaining

$$f_\delta d\alpha = i_v(df_\delta)d\alpha = i_v(df_\delta \wedge d\alpha) + df_\delta \wedge i_v d\alpha = df_\delta \wedge i_v d\alpha$$

and $(a-1)\alpha - i_v d\alpha = L_v \alpha - i_v d\alpha = di_v \alpha$.

We consider the cone $\sigma := \mathbb{R}_+ \delta$ generated by $\delta$ and define

$$\Omega_\sigma^p := \{\varphi \in \Omega_{\mathbb{C}^n}^p(\log(x_1 \ldots x_n))_0 \, | supp(\varphi) \subseteq \sigma\}$$

with the evident ($\mathbb{Q}-$) grading (compatible with $\deg f_\delta = 1$). For this grading

$$(\Omega_\sigma^p)_b = x^{bm}L(b\Delta') \otimes \Lambda^p M_{\mathbb{C}} \quad (b \in \mathbb{N}).$$

We thus have proved by Prop. 3.5:

**Corollary 3.6:** For $a \in \mathbb{N}^+$, $1 \leq r \leq n-1$, as earlier, consider $\varphi \in \Omega_0^n$ with $supp(\varphi) \subseteq a\delta^\circ$ and suppose $[\varphi] \neq 0$ in $(\Omega_\sigma^n)_a / df_\delta \wedge d(\Omega_\sigma^{n-2})_{a-1}$. Then $(\alpha, k)(\varphi) = (a-1, r-1)$.

### 3.2    Non-integral degrees

Here we consider the still missing case $a \in \mathbb{Q}\backslash\mathbb{Z}$, $1 \leq r \leq n$, and $\varphi = \varphi_1 \in \Omega_0^n$ with $supp(\varphi) \subseteq a\delta^\circ$. We show that Cor. 3.6 remains true with the appropriate changes.

Let $e := \min\{k \in \mathbb{N}^+ | kv(\mathbb{C}\{x\}) \subseteq \mathbb{Z}\}$ and $F := f(x) + y^e/e$ with the new variable $y$. In order to apply the result already proved to $F$, we consider the modified form $\phi := \varphi \wedge \frac{dy}{y} y^c$ with $c \in \mathbb{N}^+$, $a + \frac{c}{e} = [a] + 1$, which has integer degree.

The pair $(\alpha, k)$ for the original $\varphi$ is then obtained by a theorem of A. Varchenko, which is more generally concerned with a direct sum $f(x) + g(y) : (\mathbb{C}^n \times \mathbb{C}^m, 0) \to (\mathbb{C}, 0)$ of two functions with isolated singularity: Consider $\omega \in \Omega_{\mathbb{C}^n, 0}^n$, $\eta \in \Omega_{\mathbb{C}^m, 0}^m$ with their expansions

$$s[\omega] = \sum_{\alpha,k} A^{\omega}_{\alpha,k} t^{\alpha} (\ln t)^k / k!, \; s[\eta] = \sum_{\beta,l} A^{\eta}_{\beta,l} t^{\beta} (\ln t)^l / l!.$$

**Lemma 3.7:** The vector space $E_{f+g} = \Gamma(U, u^* H^{n+m-1}_{f+g})$ (cf. p. 1) for $f+g$ is isomorphic to the tensor product of the corresponding vector spaces for $f$ and $g$, and the expansion for $\omega \wedge \eta$ is

$$s[\omega \wedge \eta] = \sum_{\alpha,k} \sum_{\beta,l} \frac{1}{k!\, l!} \frac{\partial^{k+l}}{\partial \alpha^k \partial \beta^l} (t^{\alpha+\beta+1} B(\alpha+1, \beta+1) A^{\omega}_{\alpha,k} \otimes A^{\eta}_{\beta,l}$$

with $A^{\omega}_{\alpha,k} \otimes A^{\eta}_{\beta,l} \in E_f \otimes E_g \cong E_{f+g}$. (Here, $B$ is the beta-function.)

(Cf. [AGV, 13.3.5].)

Applying this to $g(y) = y^e/e$, $\omega = \varphi$, $\eta = y^c dy/y$ with

$$s[\omega] = t^{a-1} A^{\omega}_{a-1,r-1} \frac{(\ln t)^{r-1}}{(r-1)!} + \cdots, \; s[\eta] = t^{\frac{c}{e}-1} A^{\eta}_{c/e-1,0} + \cdots$$

we get

$$s[\omega \wedge \eta] = \frac{B(a,c/e)}{(r-1)!} t^{a+\frac{c}{e}-1} (\ln t)^{r-1} A^{\omega}_{a-1,r-1} \otimes A^{\eta}_{c/e-1,0} + \cdots,$$

Under the assumption for $\phi = \omega \wedge \eta$ in Cor. 3.6 the conclusion is

$$A^{\omega}_{a-1,r-1} \otimes A^{\eta}_{c/e-1,0} \neq 0,$$

in particular $(\alpha, k)(\varphi) = (a-1, r-1)$.

For the final result, we want to reformulate this assumption into one about $\varphi$.

To the face $\delta$ and the cone $\sigma$ correspond for $\Gamma_+(F) \subseteq \mathbb{R}^{n+1}$ the face $\hat{\delta} := conv(\delta \cup \{(0,e)\})$ and the cone $\hat{\sigma} := \mathbb{R}_+ \hat{\delta}$. We consider $\Omega^{n+1}_{\hat{\sigma}}/dF_{\hat{\delta}} \wedge d\Omega^{n-1}_{\hat{\sigma}}$ with the obvious grading, compatible with $\deg F_{\hat{\delta}} = 1$. Since $F_{\hat{\delta}} = f_{\delta} + y^e/e$ is invariant, the group $\mathbb{Z}/\mathbb{Z}e$ acts on each homogeneous component $(\Omega^{n+1}_{\hat{\sigma}}/dF_{\hat{\delta}} \wedge d\Omega^{n-1}_{\hat{\sigma}})_b$ by multiplication of the variable $y$ by $\zeta := e^{2\pi i/e}$ and decomposes it into the eigenspaces $(\Omega^{n+1}_{\hat{\sigma}}/dF_{\hat{\delta}} \wedge d\Omega^{n-1}_{\hat{\sigma}})_{b,c}$, with eigenvalue $\zeta^c$, $c = 0, \ldots, e-1$, of the operator $\zeta$. The $\zeta^c$−eigenspace is generated by homogeneous elements with a degree $\equiv c \bmod e$ in the variable $y$.

**Lemma 3.8:** For $a \notin \mathbb{Z}$, $a = [a] + c/e$, $b = [a] + 1$, the map

$$(\Omega^n_{\sigma}/df_{\delta} \wedge d\Omega^{n-2}_{\sigma})_a \to (\Omega^{n+1}_{\hat{\sigma}}/dF_{\hat{\delta}} \wedge d\Omega^{n-1}_{\hat{\sigma}})_{b,c}, \; [\varphi] \mapsto \left[ y^c \frac{dy}{y} \wedge \varphi \right],$$

is bijective.

Proof: We shall write $f, F$ instead of $f_{\delta}, F_{\hat{\delta}}$.

(1) The map is well-defined:

For $\varphi \in \Omega_\sigma^{n-2}$ of degree $a - 1$ holds $d\left(y^c \frac{dy}{y} \wedge \varphi\right) = -y^c \frac{dy}{y} \wedge d\varphi$ and

$$y^c \frac{dy}{y} \wedge df \wedge d\varphi = (df + y^e \frac{dy}{y}) \wedge d(y^c \frac{dy}{y} \wedge \varphi).$$

(2) The map is surjective:

The right-hand side is generated by the residue classes of $y^{c+se} \frac{dy}{y} \wedge \varphi$, $\varphi \in (\Omega_\sigma^n)_q$, $q + \frac{c}{e} + s = b$, $s \in \mathbb{N}$. As $q \neq 0$, the Poincaré-Lemma is valid, and there is a $\psi \in (\Omega_\sigma^{n-1})_q$ with $d\psi = \varphi$. For $s > 0$, the form

$$y^{c+se} \frac{dy}{y} \wedge \varphi - \left(df + y^e \frac{dy}{y}\right) \wedge d(y^{c+(s-1)e}\psi) = (c + (s-1)e)y^{c+(s-1)e} \frac{dy}{y} \wedge df \wedge \psi)$$

by now has only degree $c + (s-1)e$ in $y$, and we argue by induction. (Note that $df \wedge \varphi = 0$ by degree and a similar remark is used below.)

(3) The map is injective:

We show this by induction on $b$, beginning at $b = 1$. Let $\varphi \in (\Omega_\sigma^n)_a$, $y^c \frac{dy}{y} \wedge \varphi = dF \wedge d\Psi$ with $\Psi \in (\Omega_\sigma^{n-1})_0$. The right-hand side of this equation is invariant under the action of $\zeta$, the left-hand side belongs to eigenvalue $\zeta^c \neq 1$. So both sides and also $\varphi$ are zero.

Induction step $b - 1 \to b$: Let $\varphi \in (\Omega_\sigma^n)_a$, $\Psi \in (\Omega_\sigma^{n-1})_{b-1}$ and

(*) $\qquad y^c \frac{dy}{y} \wedge \varphi = dF \wedge d\Psi.$

We decompose $\Psi$ into homogeneous components with respect to $y$:

$$\Psi = \sum_{s \geq 0} y^{c+se} \Psi_s, \quad \Psi_s = C_s \Psi_s' + \frac{dy}{y} \wedge C_{s-1} \Psi_s'', \quad C_s := \prod_{i=0}^s (c + ie)^{-1},$$

$C_{-1} := 1$, where $\Psi_s', \Psi_s''$ contain only $dx_i/x_i$. Then

$$d\Psi = \sum_{s \geq 0} y^{c+se} (\frac{dy}{y} \wedge C_{s-1}(\Psi_s' - d\Psi_s'') + C_s d\Psi_s').$$

We deduce from (*):

$$-\varphi = df \wedge (\Psi_0' - d\Psi_0'')$$

$$d\Psi_{s-1}' = df \wedge (\Psi_s' - d\Psi_s''), s \geq 1.$$

The second equation shows, that

$$-y^c \frac{dy}{y} \wedge d\Psi_0' = dF \wedge d\widetilde{\Psi}$$

with

$$\widetilde{\Psi} = \sum_{s\geq 1} y^{c+(s-1)e}(C_{s-1}\Psi'_s + \frac{dy}{y} \wedge C_{s-2}\Psi''_s),$$

$$d\widetilde{\Psi} = \sum_{s\geq 1} y^{c+(s-1)e}(\frac{dy}{y} \wedge C_{s-2}(\Psi'_s - d\Psi''_s) + C_{s-1}d\,\Psi'_s)$$

By induction hypothesis there is a $\psi \in (\Omega_\sigma^{n-2})_{a-2}$ with

$$d\Psi'_0 = df \wedge d\psi.$$

Then $d(\Psi'_0 + df \wedge \psi) = 0$, and by the Poincaré-lemma (valid in degree $\neq 0$) there is a $\alpha \in (\Omega_\sigma^{n-2})_{a-1}$ with $\Psi'_0 + df \wedge \psi = d\alpha$, and we get $\varphi = df \wedge (-\Psi'_0 + d\Psi''_0) = df \wedge (-d\alpha + d\Psi''_0)$, as required.

With this lemma, we can express the two cases of our result in a unified form:

**Theorem 3.9:** Let $(a \in \mathbb{N}^+, 1 \leq r \leq n-1)$ or $(a \in \mathbb{Q}^+\backslash\mathbb{Z}, 1 \leq r \leq n)$. We suppose $\varphi \in \Omega_0^n$, $supp(\varphi) \subseteq a\delta^\circ$ and $[\varphi] \neq 0$ in $(\Omega_\sigma^n)_a/df_\delta \wedge d(\Omega_\sigma^{n-2})_{a-1}$. Then $(\alpha,k)(\varphi) = (a-1, r-1)$.

**Appendix 1: Complexes of differential forms**

On the open set $U = \{z | |z_i| < R\} \subseteq \mathbb{C}^n$ let $D_i$ be the divisor $(z_i)$ and for $0 \leq k \leq n$ let $D := D_1 \cup \ldots \cup D_k$ with the inclusion $j: U\backslash D \hookrightarrow U$ of the complement.

By means of the restriction maps from $U$ to $D_i$, one defines

$$\Omega^\bullet_{(U,D)} := \ker(\Omega^\bullet_U \to \oplus_{i=1}^k \Omega^\bullet_{D_i}),$$

$$\mathcal{E}^\bullet_{(U,D)} := \ker(\mathcal{E}^\bullet_U \to \oplus_{i=1}^k \mathcal{E}^\bullet_{D_i}).$$

A $p$–form $\varphi = \sum_{I,J} \varphi_{IJ}\, dz_I \wedge d\bar{z}_J$ belongs to $\mathcal{E}^\bullet_{(U,D)}$ if and only if $\varphi_{IJ}|D_i = 0$ for all $i \in [1,k]$ and $I,J$ with $i \notin I \cup J$. We have the decomposition by type

$$\mathcal{E}^p_{(U,D)} = \oplus_{r+s=p} \mathcal{E}^{rs}_{(U,D)} \text{ with } \mathcal{E}^{rs}_{(U,D)} = \mathcal{E}^{rs}_U \cap \mathcal{E}^p_{(U,D)}.$$

**Lemma A.1.1:** The exterior derivative gives rise to exact sequences

$$0 \to \Gamma j_!\mathbb{C}_{U\backslash D} \to \Gamma\Omega^\bullet_{(U,D)} \text{ and } 0 \to \Gamma j_!\mathbb{C}_{U\backslash D} \to \Gamma\mathcal{E}^\bullet_{(U,D)}.$$

Proof: (For $\mathcal{E}$, for $\Omega$ it is the same.) Starting with the ordinary Poincaré-lemma ($k = 0$) we achieve a proof by induction on $k$ if we put $D' := D_2 \cup \ldots \cup D_k$, $\widetilde{D} := D_1 \cap D'$ and construct with the inclusions $j': U\backslash D' \hookrightarrow U, \tilde{j}: D_1\backslash\widetilde{D} \hookrightarrow D_1$ the diagram

$$\begin{array}{ccccccc}
0 \to & \Gamma\mathcal{E}^\bullet_{(U,D)} & \to & \Gamma\mathcal{E}^\bullet_{(U,D')} & \to & \Gamma\mathcal{E}^\bullet_{(D_1,\widetilde{D})} & \to 0 \\
& \uparrow & & \uparrow & & \uparrow & \\
0 \to & \Gamma j_!\mathbb{C}_{U\backslash D} & \to & \Gamma j'_!\mathbb{C}_{U\backslash D'} & \to & \Gamma \tilde{j}_!\mathbb{C}_{D_1\backslash\widetilde{D}} & \to 0
\end{array}$$

which has exact lines and two exact columns by induction hypothesis.

(The terms in the lower line are the vector spaces of locally constant functions on $U\backslash D$ resp. $U\backslash D'$ resp. $D_1\backslash \widetilde{D}$, whose support is closed in $U$ resp. $U$ resp. $D_1$. They are 0 or $\mathbb{C}$, depending on whether $D$ resp. $D'$ resp. $\widetilde{D}$ is empty or not.)

**Lemma A.1.2:** The sequence

$$0 \to \Omega^p_{(U,D)} \to \mathcal{E}^{p,0}_{(U,D)} \xrightarrow{\bar{\partial}} \mathcal{E}^{p,1}_{(U,D)} \to \cdots$$

is exact.

Proof: This follows from the Dolbeaut-lemma ($k = 0$) by induction on $k$ by means of the short exact sequence

$$\begin{array}{ccccccccc}
0 & \to & \Omega^p_{(U,D)} & \to & \Omega^p_{(U,D')} & \to & \Omega^p_{(D_1,\widetilde{D})} & \to & 0 \\
& & \downarrow & & \downarrow & & \downarrow & & \\
0 & \to & \mathcal{E}^{p,\bullet}_{(U,D)} & \to & \mathcal{E}^{p,\bullet}_{(U,D')} & \to & \mathcal{E}^{p,\bullet}_{(D_1,\widetilde{D})} & \to & 0
\end{array}$$

where $D' = D_2 \cup \ldots \cup D_k$, $\widetilde{D} = D_1 \cap D'$.

For two divisors $D = D_1 \cup \ldots \cup D_k$, $Z = D_l \cup \ldots D_n$ ($0 \le k < l \le n$) we define:

$$\mathcal{E}^p_{(U,D)}(*Z) := \mathcal{E}^p_{(U,D)} \otimes \mathcal{O}(*Z),$$

$$\mathcal{E}^p_{(U,D)}(\log Z) := \mathcal{E}^p_{(U,D)}(*Z) \cap \mathcal{E}^p_U(\log Z).$$

We consider the compact subset $K := \{z \mid |z_i| \le \rho\} \subseteq U$.

**Lemma A.1.3:** Let $\varphi \in \Gamma(K, \mathcal{E}^p_{(U,D)}(*Z))$ with $d\varphi \in \Gamma(K, \Omega^{p+1}_{(U,D)}(*Z))$.

(1) $\exists\, \psi \in \Gamma(K, \mathcal{E}^{p-1}_{(U,D)}(*Z))$ with $\varphi + d\psi \in \Gamma(K, \Omega^p_{(U,D)}(*Z))$.
(2) $\exists\, \xi \in \Gamma(K, \Omega^p_{(U,D)}(*Z))$ with $d\xi = d\varphi$.

Proof: In (2) choose $\xi = \varphi + d\psi$.

(1): We argue as in [Ba3, lemme 1]. Consider $\varphi = \alpha/f$ with $\alpha \in \Gamma(K, \mathcal{E}^p_{(U,D)})$ and $f = z_l^{r_l} \ldots z_n^{r_n}$, and look at the decomposition by type $\alpha = \alpha_{i,p-i} + \cdots + \alpha_{p,0}$ ($i \ge 0$). From

$$d\varphi = \frac{d\alpha}{f} - \frac{df \wedge \alpha}{f^2} \in \Gamma(K, \Omega^{p+1}_{(U,D)}(*Z))$$

we see $\bar{\partial}\alpha_{i,p-i} = 0$. If $i < p$, by A.1.2, there exists $\beta \in \Gamma(K, \mathcal{E}^{i,p-i-1}_{(U,D)})$ with $\bar{\partial}\beta = \alpha_{i,p-i}$. Then

$$\varphi - d\frac{\beta}{f} = \frac{\alpha - d\beta}{f} + \frac{df \wedge \beta}{f^2}$$

is of type $(i+1, p-i-1) + (i+2, p-i-2) + \cdots$, and the assertion follows by induction.

**Corollary A.1.4:**

$$\Gamma(K, \Omega^\bullet_{(U,D)}(*Z)) \xrightarrow{\approx} \Gamma(K, \mathcal{E}^\bullet_{(U,D)}(*Z)).$$

**Lemma A.1.5:** Let $\varphi \in \Gamma(K, \mathcal{E}^p_{(U,D)}(*Z))$ with $d\varphi \in \Gamma(K, \mathcal{E}^{p+1}_{(U,D)}(\log Z))$.

(1) $\exists\, \psi \in \Gamma(K, \mathcal{E}^{p-1}_{(U,D)}(*Z))$ with $\varphi + d\psi \in \Gamma(K, \mathcal{E}^p_{(U,D)}(\log Z))$.
(2) $\exists\, \xi \in \Gamma(K, \mathcal{E}^p_{(U,D)}(\log Z))$ with $d\xi = d\varphi$.

The corresponding result is true for holomorphic forms.

Proof: (2) is clear by the choice $\xi = \varphi + d\psi$.

(1): Consider $\varphi = \alpha/f$ with $\alpha \in \Gamma(K, \mathcal{E}^p_{(U,D)})$ and $f = z_l^{r_l} \ldots z_n^{r_n}$. In the case $r_l, \ldots, r_n \leq 1$ there is nothing to show, since a semi-meromorphic form is logarithmic if and only if the form itself and its derivative are expressible with denominator $z_l \ldots z_n$ ([Gr-S, p. 73]). If $r_l \geq 2$ and $r_{l+1}, \ldots, r_n \geq 1$ (without loss of generality) we construct a form $\tilde\psi \in \Gamma(K, \mathcal{E}^{p-1}_{(U,D)}(*Z))$ such that in $\varphi + d\tilde\psi$ appears only denominator $z_l^{-1} f$. The assertion follows then inductively. We have $d\varphi = \frac{d\alpha}{f} - \frac{df}{f} \wedge \frac{\alpha}{f}$, and by hypothesis at least $fd\varphi \in \Gamma(K, \mathcal{E}^{p+1}_U)$, i.e.

$$\frac{df}{f} \wedge \alpha = \sum_{i=l}^n r_i \frac{dz_i}{z_i} \wedge \alpha \in \Gamma(K, \mathcal{E}^{p+1}_U).$$

By A.1.8 we can conclude from this $r_i \frac{dz_i}{z_i} \wedge \alpha \in \Gamma(K, \mathcal{E}^{p+1}_U)\ \forall\ i = l, \ldots, n$. For the vector fields $w = \frac{\partial}{\partial z_l}$, $v = z_l w$ and the inner contractions $i_w, i_v$ we have the relations

$$di_v \alpha = dz_l \wedge i_w \alpha + z_l di_w \alpha = -i_w(dz_l \wedge \alpha) + \alpha + z_l di_w \alpha,$$

$$\frac{df}{f} \wedge i_v \alpha = -i_v\left(\frac{df}{f} \wedge \alpha\right) + r_l \alpha \quad (\text{since } i_v \frac{df}{f} = r_l),$$

so

$$di_v \varphi = d\left(\frac{i_v \alpha}{f}\right) = \frac{1}{f}\left(di_v \alpha - \frac{df}{f} \wedge i_v \alpha\right) = \frac{1}{f}((1 - r_l)\alpha + z_l \beta),$$

with $\beta := di_w \alpha + i_w\left(\frac{df}{f} \wedge \alpha\right) - i_w(\frac{dz_l}{z_l} \wedge \alpha) \in \Gamma(K, \mathcal{E}^p_{(U,D)})$. Therefore, $\tilde\psi := (r_l - 1)^{-1} i_v \varphi$ is suitable. This proof also works for holomorphic forms.

**Corollary A.1.6:**

$$\Gamma(K, \Omega^\bullet_{(U,D)}(\log Z)) \xrightarrow{\approx} \Gamma(K, \Omega^\bullet_{(U,D)}(*Z)),$$

$$\Gamma(K, \mathcal{E}^\bullet_{(U,D)}(\log Z)) \xrightarrow{\approx} \Gamma(K, \mathcal{E}^\bullet_{(U,D)}(*Z)).$$

**A.1.7:** For a smooth function $f \in C^\infty(U)$ we have the equivalence:

$$f/z_1 \in C^\infty(U) \Leftrightarrow \forall\, k \geq 0: \overline{D}^k f(0, z_2, \ldots, z_n) = 0, \text{ where } \overline{D} = \partial/\partial\bar{z}_1.$$

Proof:

$\Rightarrow$: $f = z_1 g$ implies $\overline{D}^k f = z_1 \overline{D}^k g$.

$\Leftarrow$: By Taylor's formula, for any $N$ there is a representation $f = f_1 + f_2$ with

$$f_1 = \sum_{i+j \leq N} a_{ij}(z') z_1^i \bar{z}_1^j, \quad f_2 = \sum_{i+j=N+1} a_{ij}(z) z_1^i \bar{z}_1^j$$

In the first sum $z' = (z_2, \ldots, z_n)$ and $i \geq 1$ by assumption. The claim follows from $\bar{z}_1^{N+1}/z_1 \in C^{N-1}(U)$.

**A.1.8:** Assume a relation $\sum_{i=1}^n \frac{f_i}{z_i} = g$ for smooth functions $f_1, \ldots, f_n, g \in C^\infty(U)$. Then $f_i/z_i \in C^\infty(U)\ \forall\, i = 1, \ldots, n$.

Proof: With $p_i := z_1 \ldots \hat{z}_i \ldots z_n$, $p = z_1 \ldots z_n$ the equality $p_i \overline{D}_i^k f_i = -\sum_{j \neq i} p_j \overline{D}_i^k f_j + p \overline{D}_i^k g$ shows that $\overline{D}_i^k f_i = 0$ for $z_i = 0$.

A variant of A.1.6 is the following:

**Lemma A.1.9:**

$$\Gamma(K, \Omega_U^\bullet(\log D \cup Z)) \xrightarrow{\approx} \Gamma(K, \Omega_U^\bullet(\log D)(* Z)),$$

$$\Gamma(K, \mathcal{E}_U^\bullet(\log D \cup Z)) \xrightarrow{\approx} \Gamma(K, \mathcal{E}_U^\bullet(\log D)(* Z)).$$

Proof: In the smooth case as a template. Consider $\varphi \in \Gamma(K, \mathcal{E}_U^p(\log D)(* Z))$ with $d\varphi \in \Gamma(K, \mathcal{E}_U^{p+1}(\log D \cup Z))$. We modify the proof A.1.5, (1), to find a $\psi \in \Gamma(K, \mathcal{E}_U^{p-1}(\log D)(* Z))$ with $\varphi + d\psi \in \Gamma(K, \mathcal{E}_U^p(\log D \cup Z))$. Let $\varphi = \alpha/f$ with $\alpha \in \Gamma(K, \mathcal{E}_U^p(\log D))$, $f = z_l^{r_l} \ldots z_n^{r_n}$, $r_{l+1}, \ldots, r_n \geq 1$, $r_l \geq 2$ (otherwise $\varphi$ would be logarithmic). We can see as before

$$\sum_{i=l}^n r_i \frac{dz_i}{z_i} \wedge \alpha \in \Gamma(K, \mathcal{E}_U^{p+1}(\log D)) \text{ and } r_i \frac{dz_i}{z_i} \wedge \alpha \in \Gamma(K, \mathcal{E}_U^{p+1}(\log D))\ (1 \leq i \leq n).$$ Then

$$\beta = d i_w \alpha + i_w\left(\frac{df}{f} \wedge \alpha\right) - i_w\left(\frac{dz_l}{z_l} \wedge \alpha\right) \in \Gamma(K, \mathcal{E}_U^p(\log D)),$$

and the rest of the proof remains the same.

We have a similar result with arbitrary singularities along $Z$. Let $j: U \backslash Z \hookrightarrow U$ be the inclusion.

**Lemma A.1.10:**

$$\Gamma(K, \Omega^\bullet_{(U,D)}(\log Z)) \xrightarrow{\approx} \Gamma(K, j_*\Omega^\bullet_{(U\setminus Z, D\setminus Z)})$$
$$\downarrow \approx \qquad\qquad \downarrow \approx$$
$$\Gamma(K, \mathcal{E}^\bullet_{(U,D)}(\log Z)) \xrightarrow{\approx} \Gamma(K, j_*\mathcal{E}^\bullet_{(U\setminus Z, D\setminus Z)})$$

Proof: The left column is quasi-isomorphic by A.1.4 and A.1.6. The right column is quasi-isomorphic by A.1.1 and the fact that $U\setminus Z$ is a Stein-space. In the first line we can reduce to the case $D = \emptyset$ by the exact sequence

$$0 \to \Omega^\bullet_{(U,D)}(\log Z) \to \Omega^\bullet_{(U,D')}(\log Z) \to \Omega^\bullet_{(D_1,\tilde{D})}(\log \tilde{Z}) \to 0$$

($D' = D_2 \cup \ldots \cup D_k$, $\tilde{D} = D' \cap D_1$, $\tilde{Z} = Z \cap D_1$) and applying $\Gamma(K, .)$ and $\Gamma(K, j_*.)$. In this case it follows by expansion into Laurent-series, that the $p$-th cohomology of the two terms has the basis $(2\pi i)^{-p}(dz)_I/z_I$, $I \subseteq [l, N]$, $|I| = p$.

**Appendix 2: Hypersurfaces in torus embeddings**

Let $M \cong \mathbb{Z}^n$, $n \geq 1$, with the dual lattice $N := Hom(M, \mathbb{Z})$, $T := Spec\ \mathbb{C}[M]$, $g \in \mathbb{C}[M]$, and $\Delta := \Delta(g) = conv(supp(g)) \subseteq M_\mathbb{R} := M \otimes \mathbb{R}$ the Newton polyhedron.

Let $\Sigma$ be a fan on $N_\mathbb{R} := N \otimes \mathbb{R}$ such that the supporting function $s_\Delta(a) := min\{\langle a, x\rangle | x \in \Delta\}$ is linear on all $\sigma \in \Sigma$ (in case $\dim \Delta = n$, there is a coarsest such fan $\Sigma_\Delta$). Let $X = X_\Sigma$ be the corresponding compact torus embedding (in case $\Sigma = \Sigma_\Delta$ it is projective).

Then the closure $\bar{V}$ of the divisor $V = (g)$ on $T \subseteq X$ is again a hypersurface (i.e. effective Cartier-divisor), which contains no stratum of $X$.

The intersections of the strata of $X$ with $\bar{V}$ are smooth hypersurfaces within these, if and only if $g$ is globally non-degenerate, i.e. if for all faces $\delta$ of $\Delta$ ($\Delta$ included) the Laurent-polynomial

$$g_\delta := \sum_{m \in \delta} g_m\, x^m, \text{ where } g = \sum_{m \in M} g_m\, x^m,$$

defines on $T$ a smooth hypersurface. In this case we call $Z = \bar{V} \subseteq X$ a non-degenerate hypersurface.

Now we assume that $\dim \Delta = n$, $\Sigma$ is regular, i.e. $X$ is smooth, and $Z = \bar{V} \subseteq X$ non-degenerate. Then $D \cup Z$, $D := X\setminus T$, is a divisor with normal crossings. We decompose: $(g) = Z - C$, $|C| \subseteq D$. The sheaf $\mathcal{O}(C)$ is a $T$-invariant (invertible) fractional ideal with order-function $ord\ \mathcal{O}(C) = s_\Delta$.

Let $j: X\setminus(D \cup Z) \hookrightarrow X$, $j_0: X\setminus Z \hookrightarrow X$ be the inclusions and $c_X := \{A | A \subseteq X \text{ compact}\}$, $\phi := c_X \cap (X\setminus Z)|T\setminus V$ supporting families.

**Lemma A.2.1:** For all $p \geq 0$ we have:

$$H^p(T\backslash V) = \mathbb{H}^p(\Omega_X^\bullet(\log D)(*Z)) = H^p \, \Gamma(X, \Omega_X^\bullet(\log D)(*Z)),$$

$$H_\phi^p(T\backslash V) = \mathbb{H}^p(\Omega_{(X,D)}^\bullet(*Z)) = H^p \, \Gamma(X, \Omega_{(X,D)}^\bullet(*Z)).$$

Proof: The left equality follows as in section 2.2.2 from

$$\Omega_X^\bullet(\log D)(*Z) \approx \Omega_X^\bullet(\log D \cup Z) \approx \mathcal{E}_X^\bullet(\log D \cup Z) \approx j_*\mathcal{E}_{X\backslash(D\cup Z)}^\bullet$$

(A.1.9, A.1.10), resp. $\Omega_{(X,D)}^\bullet(*Z) \approx j_{0*}j_0^{-1}\mathcal{E}_{(X,D)}^\bullet$ (A.1.6, A.1.10). It suffices then to show, that $\Omega_X^\bullet(\log D)(*Z)$ and $\Omega_{(X,D)}^\bullet(*Z)$ are acyclic. If we apply the formula

$$H^i(X, \mathcal{F}(*Z)) = \varinjlim_{l \geq 0} H^i(X, \mathcal{F}(lZ))$$

for coherent sheaves on a compact space (cf. [Go, p. 194]), this follows from:

**Lemma A.2.2:** ([TE], [Da 1]) Let $n \geq 1$.

a) $H^0(X, \mathcal{O}(lC)) = L(l\Delta), l \geq 0$;
   $H^n(X, \mathcal{O}(lC)) \cong L(l\Delta^\circ), l < 0$;
   $H^i(X, \mathcal{O}(lC)) = 0$, in the remaining cases;
b) $H^0(X, \mathcal{O}(lC - D)) = L(l\Delta^\circ), l > 0$;
   $H^n(X, \mathcal{O}(lC - D)) \cong L(l\Delta), l \leq 0$;
   $H^i(X, \mathcal{O}(lC - D)) = 0$, in the remaining cases;
   (with the notation $L(A) := \{g \in \mathbb{C}[M] | supp(g) \subseteq A\}$ for $A \subseteq \mathbb{R}^n$).
c) $\Omega_X^p(\log D) \cong \mathcal{O}_X \otimes_\mathbb{Z} \Lambda^p M$, in particular: $\omega_X = \Omega_X^n = \mathcal{O}(-D) \otimes_\mathbb{Z} \Lambda^n M$;
   $\Omega_{(X,D)}^p = \Omega_X^p(\log D)(-D) \cong \omega_X \otimes_\mathbb{Z} \Lambda^p M$.

Proof: Let $h = s_\Delta: |\Sigma| \to \mathbb{R}$ be the order function of $\mathcal{O}(C)$.

a): By [TE], [Da 1] the $m$-th ($m \in M$) component of $H^i(X, \mathcal{O}(lC))$ can be calculated as

$$H^i(X, \mathcal{O}(lC))(m) \cong H^i_{Z_m}|\Sigma| \text{ with } Z_m := \{a \in |\Sigma| \,|\, m(a) \geq lh(a)\}$$

More explicitly,

$$H^0_{Z_m}|\Sigma| = \ker(H^0|\Sigma| \to H^0(|\Sigma|\backslash Z_m) \cong \begin{cases} \mathbb{C}, Z_m = |\Sigma| \\ 0, \text{ else} \end{cases}$$

$$H^i_{Z_m}|\Sigma| = \tilde{H}^{i-1}(|\Sigma|\backslash Z_m), i \geq 1, \text{ (reduced cohomology)}.$$

In the case $l < 0$, $Z_m$ is a convex cone with vertex, and

$$|\Sigma|\backslash Z_m \simeq \begin{cases} \text{point}, Z_m \neq \{0\} \Leftrightarrow -m \notin -l\Delta^\circ \Leftrightarrow m \notin l\Delta^\circ \\ S^{n-1}, Z_m = \{0\} \Leftrightarrow m \in l\Delta^\circ \end{cases}.$$

The assertions follow from this.

b) A direct proof without Serre-duality is like this: For $m \in M$ we put

$$A := A_m := \{a \in |\Sigma| \, | \, m(a) \geq ord \, \mathcal{O}(lC - D)(a)\}$$

$$B := B_m := \{a \in |\Sigma| \, | \, m(a) > lh(a)\} \cup \{0\}.$$

Then

$$H^i(X, \mathcal{O}(lC - D))(m) \cong H^i_{A_m}|\Sigma|.$$

We have for $\sigma \in \Sigma$: $H^0(\sigma \backslash A) = H^0(\sigma \backslash B)$, $H^i(\sigma \backslash A) = H^i(\sigma \backslash B) = 0$, $i > 0$, because of

$$\sigma \subseteq A \Leftrightarrow m \geq lh + ord \, \mathcal{O}(-D) \text{ on } \sigma \Leftrightarrow m > lh \text{ on } \sigma \backslash \{0\} \Leftrightarrow \sigma \subseteq B$$

(where we have used, that $ord \, \mathcal{O}(-D)(a) = 1$ for all primitive edge vectors $a$ of $\sigma$).

$H^i(|\Sigma|\backslash A)$ and $H^i(|\Sigma|\backslash B)$ are the Cech-cohomology of the acyclic coverings $(\sigma\backslash A)_{\sigma \in \Sigma}$ resp. $(\sigma\backslash B)_{\sigma \in \Sigma}$ ([Go, p. 209]) and are thus equal. Then also $H^i_A|\Sigma|$ and $H^i_B|\Sigma|$ are equal, and we may compute $H^i_B|\Sigma|$.

In the case $l > 0$ the set $|\Sigma|\backslash B_m$ is a convex cone with vertex without $\{0\}$, and

$$H^0_{B_m}|\Sigma| = \ker(H^0|\Sigma| \to H^0(|\Sigma|\backslash B_m)) \cong \begin{cases} \mathbb{C}, & B_m = |\Sigma| \Leftrightarrow m \in l\Delta^\circ \\ 0 & else \end{cases}$$

$$H^i_{B_m}|\Sigma| = \tilde{H}^{i-1}(|\Sigma|\backslash B_m) = 0, i \geq 1.$$

In the case $l \leq 0$, $B_m$ is an open convex cone with vertex $\cup \{0\}$ or an open half space $\cup \{0\}$, and we have $B_m = \{0\} \Leftrightarrow -m \in -l\Delta \Leftrightarrow m \in l\Delta$. This shows

$$H^0_{B_m}|\Sigma| = 0, \quad H^i_{B_m}|\Sigma| \cong \begin{cases} \mathbb{C}, m \in l\Delta, \ i = n \\ 0 \ else \end{cases} \text{ for } i \geq 1,$$

as required.

c) The formula about $\Omega^p_X(\log D)$ is easy to verify on toric charts of $X$, cf. [Da 1], $\Omega^p_{(X,D)} = \Omega^p_X(\log D)(-D)$ is true by definition.

**Lemma A.2.3:** The exterior derivative $d: \Omega^p_X(\log D)(lZ) \to \Omega^{p+1}_X(\log D)((l+1)Z)$ induces $\mathcal{O}_X$-linear exact sequences

$$0 \to \Omega^0_X(\log D) \otimes \mathcal{N} \to \cdots \to \Omega^n_X(\log D) \otimes \mathcal{N}^{n+1} \to 0,$$

$$0 \to \Omega^0_{(X,D)} \otimes \mathcal{N} \to \cdots \to \Omega^n_{(X,D)} \otimes \mathcal{N}^{n+1} \to 0,$$

where $\mathcal{N} = \mathcal{O}(Z) \otimes \mathcal{O}_Z$ is the normal bundle.

Proof: Let $U \subseteq X$ be open with a coordinate system $z, z_2, \ldots, z_n$ such that $D \cup Z|U$ is a union of coordinate hypersurfaces and $Z|U = (z)$. Then we have (for $\alpha$ out of $\Omega^{l-1}_U(\log D), l \geq 1$)

$$d\left(\frac{\alpha}{z^l}\right) = \frac{1}{z^l}\left(-l\frac{dz}{z}\wedge \alpha + d\alpha\right) \equiv -l\frac{dz}{z}\wedge \frac{\alpha}{z^l} \mod \left(\frac{1}{z^l}\Omega^l(\log D)\right).$$

Furthermore, $w = \frac{\partial}{\partial z}$ is a vector field with $i_w\Omega^p_U(\log D) \subseteq \Omega^{p-1}_U(\log D)$, and for $v = zw$ we have

$$-\left(\frac{1}{l}i_v d + \frac{1}{l-1}di_v\right)\left(\frac{\alpha}{z^l}\right) \equiv i_v\left(\frac{-l}{-l}\frac{dz}{z}\wedge\frac{\alpha}{z^l}\right) + \left(\frac{-(l-1)}{-(l-1)}\frac{dz}{z}\wedge\frac{i_w\alpha}{z^{l-1}}\right) \mod \left(\frac{1}{z^{l-1}}\right) \equiv \frac{\alpha}{z^l},$$

where the terms with denominator $l - 1$ are to be omitted if $l = 1$. This shows exactness of the first sequence. The second sequence is the tensor product of the first one with $\mathcal{O}(-D)$ and remains exact.

Of course, also tensoring with $\mathcal{N}^r, r \in \mathbb{Z}$, or multiplying the differentials by non-zero constants does not affect exactness of the sequences.

**Lemma A.2.4:** The sequences of section modules (for $r \in \mathbb{Z}$)

$$0 \to \Gamma(X, \Omega^0_X(\log D) \otimes \mathcal{N}^{-r}) \to \cdots \to \Gamma(X, \Omega^n_X(\log D) \otimes \mathcal{N}^{-r+n}) \to 0,$$

$$0 \to \Gamma(X, \Omega^0_{(X,D)} \otimes \mathcal{N}^{1-r}) \to \cdots \to \Gamma(X, \Omega^n_{(X,D)} \otimes \mathcal{N}^{1-r+n}) \to 0,$$

are exact, except at the last term, and exact everywhere if $r < 0$.

Proof: The underlying sheaf sequence $\mathcal{S}^\bullet$ (with differential given by $\frac{dz}{z}\wedge$ for a local equation $z$ of $Z$) is a resolution of zero. If we consider the first spectral sequence of hypercohomology, we therefore have:

$$'E_1^{pq} = H^q(X, \mathcal{S}^p), \quad 'E_\infty^{pq} = 0 \text{ for all } p, q.$$

By Lemma A.2.2 (together with the sequence $0 \to \mathcal{O}((l-1)Z) \to \mathcal{O}(lZ) \to \mathcal{N}^l \to 0$ and $\dim Z = n - 1$) only the marked $'E_1$-terms may be non-zero.

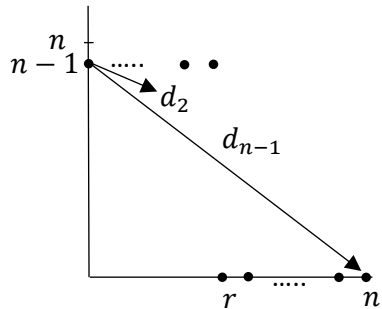

Keeping in mind that $d_s$ is of degree $(s, 1 - s)$, the diagram shows:

$$'E_2^{p,0} = 'E_\infty^{p,0} = 0 \text{ for } p \neq n, \text{ resp. for } p \text{ arbitrary if } r < 0,$$

as was required.

To simplify notation, we write $\Omega^\bullet$ for the complex $\Omega_X^\bullet(\log D)$, resp. $\Omega_{(X,D)}^\bullet$, and $K^\bullet$ for $\Omega^\bullet(*Z)$. We define a filtration of $K^\bullet$ by the subcomplexes $F^p K^\bullet$:

$$0 \to \Omega^p(Z) \to \Omega^{p+1}(2Z) \to \cdots.$$

In the diagram

$$0 \to \Gamma(F^{p+1}K^\bullet) \to \Gamma(F^p K^\bullet) \to \Gamma(Gr^p K^\bullet) \to 0,$$

the lines are exact by A.2.2, and the right column begins with

$$0 \to \Gamma(\Omega^p(Z)) \xrightarrow{\bar{d}} \Gamma(\Omega^{p+1} \otimes \mathcal{N}^2) \to \cdots.$$

Likewise, $\Gamma(\Omega^p(Z))/\Gamma(\Omega^p) \subseteq \Gamma(\Omega^p \otimes \mathcal{N})$, with equality, if $\Omega^p = \Omega_X^p(\log D)$ or $n \geq 2$.

By A.2.4 we have

$$H^q \Gamma(Gr^p K^\bullet) = 0 \quad \text{for } p < q < n \text{ or } (p < q, \text{ if } p \leq 0, \text{ resp. } p < 0).$$

In any case, the connecting operator $\delta: H^q\Gamma(Gr^p K^\bullet) \to H^{q+1}\Gamma(F^{p+1}K^\bullet)$ is zero, since (in the above sequence) $\ker \bar{d} = \Gamma(\Omega^p)$ for $p < n$, and $d\Gamma(\Omega^p) = 0$.

We have obtained:

**Lemma A.2.5:** The canonical map

$$H^q\Gamma(X, F^{p+1}\Omega^\bullet(*Z)) \to H^q\Gamma(X, F^p\Omega^\bullet(*Z)),$$

where $\Omega^\bullet = \Omega_X^\bullet(\log D)$, resp. $\Omega_{(X,D)}^\bullet$, is injective for all $p, q \in \mathbb{Z}$, and surjective for $p < q < n$ or $(q = n, p \leq 0$, resp. $p < 0)$.

**Corollary A.2.6:** For $l \in \mathbb{N}^+$ there are inclusions

$$\frac{\Gamma(X, \Omega^n(\log D)(lZ))}{d\Gamma(X, \Omega^{n-1}(\log D)((l-1)Z))} \subseteq H^n(T\backslash V) \quad \text{(with equality for } l \geq n\text{)},$$

$$\frac{\Gamma(X, \Omega_{(X,D)}^n(lZ))}{d\Gamma(X, \Omega_{(X,D)}^{n-1}((l-1)Z))} \subseteq H_\phi^n(T\backslash V) \quad \text{(with equality for } l > n\text{)}.$$

The quotients of these filtrations (in case $n \geq 2$) are just the $n$-th cohomology in A.2.4 and are easy to calculate.

Finally, we have to provide the proofs for Lemmas 3.3 and 3.4.

As before, let $i^*: H^{n-1}T \to H^{n-1}V$ the map induced by inclusion, $c: H_c^{n-1}V \to H^{n-1}V$ the canonical map, and $R_c: H_\phi^n(T\backslash V) \to H_c^{n-1}V$ the residue map. We need to show, that $i^*$ is injective.

We have the representations $H^k T = \Gamma(X, \Omega^k(\log D)) = \Lambda^k M_\mathbb{C}$ and $\Gamma(X, \Omega^n(\log D)(Z)) = g^{-1}L(\Delta) \otimes \Lambda^n M_\mathbb{C} \subseteq H^n(T\backslash V)$ (by A.2.2, A.2.6). We can construct a commutative diagram:

$$\begin{array}{ccc}
& H^{n-1}(T) & \\
dg/g \downarrow & & \searrow i^* \\
H^n(T) \to & H^n(T\backslash V) \xrightarrow{R} & H^{n-1}(V)
\end{array}$$

Since $\dim \Delta = n$, the derivatives $x_1 \frac{\partial g}{\partial x_1}, \ldots, x_n \frac{\partial g}{\partial x_n}$ and $g$ ($x_1, \ldots, x_n$ a character basis for $T$) are linear independent in $L(\Delta)$. This means, that the intersection of the images of the injective maps

$$\Lambda^n M_{\mathbb{C}} \xrightarrow{g} L(\Delta) \otimes \Lambda^n M_{\mathbb{C}} \xleftarrow{dg \wedge} \Lambda^{n-1} M_{\mathbb{C}}$$

is zero, and we have shown:

**Lemma A.2.7:** The restriction $i^*: H^{n-1}T \to H^{n-1}V$ is injective.

**Lemma A.2.8:** For $n \geq 2$ we have $image(c) \cap image(i^*) = 0$ in $H^{n-1}(V)$, for $n = 1$ we have $image(R_c) \cap image(i^*) = 0$ in $H^0(V) = H_c^0(V)$.

Proof for A.2.8: First assume $n \geq 2$. We consider the set $E$ of primitive edge vectors in $\Sigma$, the subset $E_0$ of those in $\Sigma_\Delta$, and the associated divisors $D(e)$ and embedded tori $T_e \subseteq D(e)$, $e \in E$. For $e \in E_0$, $Z \cap D(e) \subseteq D(e)$ is a non-degenerate hypersurface, whose Newton polyhedron is (essentially) the face of $\Delta$ corresponding to $e$, which is of dimension $n - 1 \geq 1$. Let $V_e = Z \cap T_e$. The torus $T_e$ is the restriction of the divisor $D(e)$ to $(X \setminus \bigcup_{\substack{e' \in E, \\ e' \neq e}} D(e'))$,

$$T_e = (X \setminus \bigcup_{\substack{e' \in E, \\ e' \neq e}} D(e')) \cap D(e),$$

and the complement is $T$. In the following diagram, we take the corresponding residue map:

$$\begin{array}{ccc} H^{n-1}(T) & \xrightarrow{(R_{T_e})} & \bigoplus_{e \in E_0} H^{n-2}(T_e) \\ i^* \downarrow & & \downarrow (i|V_e)^* \\ H_c^{n-1}(V) \xrightarrow{c} H^{n-1}(V) & \xrightarrow{(R_{V_e})} & \bigoplus_{e \in E_0} H^{n-2}(V_e) \end{array}$$

The columns are injective by A.2.7. The map $(R_{T_e})$ is described by $(i_e): \Lambda^{n-1} M_{\mathbb{C}} \to \bigoplus_{e \in E_0} \Lambda^{n-2} \ker e$ and is injective (noticing $\dim \Delta = n \geq 2$). The composition $(R_{V_e}) \circ c$ is zero, because $c$ factorizes over the restriction $H^{n-1}(Z) \to H^{n-1}(V)$.

This shows the first assertion. The second assertion for $n = 1$ follows from the fact that the sum of the residues of a meromorphic $1$-form on $X$ is zero.

**Lemma A.2.9:** The intersection of the images of the canonical maps

$$H^n(T) \to H^n(T \setminus V) \leftarrow H_\phi^n(T \setminus V)$$

is zero for $n \geq 1$.

Proof: We have a diagram of residue maps

$$H^n(T) \xrightarrow{(R_{T_e})} \bigoplus_{e \in E_0} H^{n-1}(T_e)$$
$$j^* \downarrow \qquad\qquad\qquad \downarrow (j|T_e \backslash V_e)^*$$
$$H^n_\phi(T \backslash V) \to H^n(T \backslash V) \xrightarrow{(R_{T_e}|T \backslash V)} \bigoplus_{e \in E_0} H^{n-1}(T_e \backslash V_e)$$

where the columns are injective by A.2.6. The map $(R_{T_e})$ is described by $(i_e): \Lambda^n M_{\mathbb{C}} \to \bigoplus_{e \in E_0} \Lambda^{n-1} \ker e$ and is also injective. Again the composition in the bottom line is zero, since the left map factorizes over the restriction $H^n(X \backslash Z) \to H^n(T \backslash V)$.